\ifdefined\pdftexversion\pdfoutput=1\fi
\documentclass[ijoc,nonblindrev]{informs4}

\usepackage[letterpaper,margin=1in,headsep=18pt]{geometry}
\usepackage[T1]{fontenc}
\usepackage[utf8]{inputenc}
\usepackage[varqu,scaled=0.95]{zi4} %

\makeatletter
\renewcommand\AAnormalsizeXI{%
   \@setfontsize\normalsizeXI{11}{14.5}%
   \Dskips{1.0}{0.3}{0.5}{0.6}%
   \let\@listi\@listI
   \baselineskip14.5pt plus 0.5pt minus 0.3pt
   \@bls=14.5pt\@em=11pt\relax}
\renewcommand\AAsmallX{%
   \@setfontsize\smallX{10}{12}%
   \Dskips{0.8}{0.3}{0.5}{0.6}%
   \def\@listi{\leftmargin\leftmargini
     \topsep 0\p@ \parsep 2\p@ \itemsep \parsep}%
   \baselineskip12pt plus 0.5pt minus 0.3pt
   \@bls=12pt\@em=10pt\relax}
\renewcommand\AAfootnotesizeIX{%
   \@setfontsize\footnotesizeIX\@ixpt\@xipt
   \Dskips{0.8}{0.3}{0.5}{0.6}%
   \def\@listi{\leftmargin\leftmargini
      \topsep 0\p@ \parsep 0\p@ \itemsep \parsep}%
   \baselineskip11pt plus 0.5pt minus 0.3pt
   \@bls=11pt\@em=9pt\relax}
\makeatother
\OneAndAHalfSpacedXI
\newcommand{\coderepo}{\url{https://github.com/Zhengzhong-You/LDST}}
\makeatletter
\long\def\theARTICLETOPLEFT{{\tabcolsep0pt\fs.9.11.\sffamily
  \begin{tabular}[t]{@{}l}%
  {\bfseries Preprint}\\[1pt]
  {\fs.8.10.\mdseries October 2026}
  \end{tabular}}}
\def\theARTICLETOPRIGHT{\llap{{\tabcolsep0pt\fs.8.10.\sffamily
  \begin{tabular}[t]{r@{}}%
  {\bfseries Code and data}\\[1pt]
  {\mdseries\coderepo}
  \end{tabular}}}}
\RRHFirstLine{\bf\theRUNAUTHOR:\enskip {\it\theRUNTITLE}}
\RRHSecondLine{\fs.7.9.\rm Preprint}
\LRHFirstLine{\bf\theRUNAUTHOR:\enskip {\it\theRUNTITLE}}
\LRHSecondLine{\fs.7.9.\rm Preprint}
\makeatother

\usepackage{natbib}
\bibpunct[, ]{(}{)}{,}{a}{}{,}
\def\bibfont{\small}

\usepackage{graphicx}
\usepackage{booktabs}
\usepackage{multirow}
\usepackage{array}
\usepackage{float}
\usepackage{mathtools}
\usepackage{enumitem}
\usepackage{pifont}
\usepackage{tcolorbox}
\usepackage{xcolor}
\usepackage{xurl}
\usepackage{microtype}
\usepackage[colorlinks=true,linkcolor=blue!50!black,citecolor=blue!50!black,urlcolor=blue!50!black]{hyperref}
\hypersetup{
  pdftitle={Log-Driven Semantic Tuning for Mixed Integer Programming Solvers},
  pdfauthor={Zhengzhong Ricky You, Kaizhao Sun, Zeping Min, Wotao Yin, Xinshang Wang},
  pdfkeywords={mixed-integer programming, algorithm configuration, large language models, parameter tuning, SCIP}
}

\newcolumntype{L}[1]{>{\raggedright\arraybackslash}m{#1}} %
\newcolumntype{C}[1]{>{\centering\arraybackslash}m{#1}} %

\makeatletter
\def\FigureCaptionFontStyle{\rmfamily\mdseries}
\def\TableCaptionFontStyle{\rmfamily\mdseries}
\def\FigureNameFontStyle{\sffamily\bfseries}
\def\TableNameFontStyle{\sffamily\bfseries}
\long\def\@maketablecaption#1#2{\EGT\TableCaptionFontStyle
  \setbox9=\hbox{\TableNameFontStyle #1\kern10pt}\hangindent\wd9\noindent
  \hbox{\TableNameFontStyle #1\kern10pt}\ignorespaces #2\endgraf\vskip6pt}
\long\def\@makefigurecaption#1#2{\EGT\FigureCaptionFontStyle
  \setbox9=\hbox{\FigureNameFontStyle #1\kern10pt}\hangindent\wd9\noindent
  \hbox{\FigureNameFontStyle #1\kern10pt}\ignorespaces #2\endgraf\vskip6pt}
\makeatother

\newcommand{\name}{LDST}

\newcommand{\xmark}{\ding{55}}

\TheoremsNumberedThrough
\EquationsNumberedThrough

\makeatletter
\newcommand{\appendix}{\par
  \setcounter{section}{0}\setcounter{subsection}{0}\setcounter{subsubsection}{0}%
  \gdef\thesection{\@Alph\c@section}%
  \gdef\theHsection{appendix.\@Alph\c@section}}
\makeatother

\newcommand{\floatnotes}[1]{%
  \par\smallskip
  \begin{minipage}{\linewidth}
  \TableNoteStyle\raggedright\noindent
  \NOTE{Notes.}#1\par
  \end{minipage}}

\begin{document}

\RUNAUTHOR{You et al.}
\RUNTITLE{Log-Driven Semantic Tuning for MIP Solvers}

\TITLE{Log-Driven Semantic Tuning for Mixed Integer Programming Solvers}

\ARTICLEAUTHORS{%
\AUTHOR{Zhengzhong Ricky You, Kaizhao Sun, Zeping Min, Wotao Yin, Xinshang Wang\textsuperscript{*}}
\AFF{Alibaba Group US, Bellevue, WA, USA\\
\textsuperscript{*}Correspondence: Xinshang Wang, \EMAIL{xinshang.w@alibaba-inc.com}}
}

\ABSTRACT{%
Automated configuration of mixed-integer programming (MIP) solvers largely relies on black-box optimization. Each trial requires an expensive solver run and provides only a scalar outcome, so tuning can require many such evaluations. Human experts take a different approach. They inspect solver logs, identify performance bottlenecks, and adjust the parameters that address them. We propose Log-Driven Semantic Tuning (\name{}), which uses a large language model (LLM) to automate this process. \name{} extracts compact diagnostics from solver logs, uses the LLM to propose sparse edits to the current best configuration based on parameter meanings and observed solver behavior, and validates each proposal before execution. A regime-adaptive objective selects the tuning signal according to the information available from a time-limited run, while a log-based completion signal estimates completion time from branch-and-bound tree growth recorded in ordinary solver logs, so that solved and unsolved runs share one time scale. On 215 instances from five benchmark sets, \name{} reduces geometric-mean loss by $43.9\%$ on easy instances and $16.7\%$ on intermediate instances relative to the best of \textsc{SMAC3}, \textsc{RBFOpt}, and random search. In both regimes, \name{} achieves a lower loss with 10 proposals than any of these baselines achieves with 50. Configurations that \name{} tunes on a single instance also transfer better than baseline configurations to unseen instances of the same family.}


\maketitle

\section{Introduction}

The performance of mixed-integer programming (MIP) solvers such as \textsc{SCIP} \citep{achterberg2009scip} depends heavily on parameter configuration, and instance-specific or class-specific tuning can yield order-of-magnitude speedups \citep{hutter2010automated}. Yet most automated configuration methods rely on black-box optimization (BBO; \citealt{schede2022survey,hoos2012automated,
hutter2009paramils,hutter2011sequential,lindauer2022smac3,lopez2016irace,
ansotegui2009gender,xu2011hydra,
hutter2014aclib,eggensperger2019pitfalls}). These methods typically use only scalar feedback, such as runtime or a derived loss, making the search highly sample inefficient when solver evaluations are expensive. A run that times out because of weak heuristics may be difficult to distinguish from one slowed by excessive linear programming (LP) relaxation work without additional evaluations. Human experts take a different approach. They examine solver logs, infer performance bottlenecks from process signals, and make targeted parameter changes. For example, a run that finds primal solutions early but makes little progress on the dual bound requires a different intervention from one dominated by the computational cost of strong branching.

We show that solver-process information can systematically guide algorithm configuration. We propose Log-Driven Semantic Tuning (\name{}), which extracts compact diagnostics from solver logs and uses a large language model (LLM) to propose sparse edits to the current best configuration. Each edit changes only a limited number of parameters from a predefined whitelist, while all other parameter values are inherited from the current best configuration. The tuning loop combines diagnostic summaries, memory centered on the best configuration, and explicit validation of the LLM's advice, allowing semantic reasoning to guide the search while each validated configuration is evaluated through an actual solver run. By modifying the current best configuration rather than constructing each candidate from scratch, \name{} maintains a local and stable search, makes individual parameter changes easier to interpret, and uses information from earlier runs to inform subsequent decisions.

The tuning loop also uses a regime-adaptive objective, as the most informative progress signal depends on the outcome of the default run at the time limit. Case~1 contains easy instances solved by the default configuration, Case~2 contains intermediate instances with a feasible solution but no proof of optimality, and Case~3 contains hard instances with only dual-bound progress. For easy instances, solve time is the natural objective, but a configuration that reaches the time limit has no observed solve time. Solved status, final gap, and the primal--dual integral cannot be compared with observed solve times, and existing tree-size estimators are built into specific solvers \citep{hendel2022estimating}. We therefore introduce a log-based completion signal that estimates the completion time of an unfinished run from branch-and-bound tree growth, which MIP solver logs commonly report. Solved and unsolved runs then share one time scale, so a single Case~1 objective ranks both. In a retrospective ranking study, the signal also orders unfinished runs more accurately than a proxy for the primal--dual integral. Across 215 instances from five benchmark sets, \name{} achieves final geometric-mean losses of $0.274$ and $0.557$ in Cases~1 and~2 after 50 proposals, respectively, $43.9\%$ and $16.7\%$ below the best scalar-feedback baselines ($0.490$ and $0.668$). By trial~10, \name{} already reaches losses of $0.396$ and $0.660$, lower than any baseline achieves after 50 proposals, while the corresponding baseline losses at trial~10 remain at $0.588$--$0.706$ and $0.763$--$0.807$.

We make the following contributions.
\begin{itemize}
   \item \textit{We propose \name{}, a log-driven semantic tuner for MIP solvers}. \name{} extracts compact diagnostics from solver logs and uses an LLM to propose sparse edits to the current best configuration. Each proposal undergoes schema, domain, novelty, and settings-loading checks before evaluation. This provides a process-aware alternative to scalar-feedback tuning without requiring task-specific training, solver-source modifications, or access to internal solver states.

   \item \textit{We introduce a regime-adaptive tuning objective with a log-based completion signal for easy instances}. The objective uses estimated completion time for easy instances, the primal--dual integral \citep{berthold2013measuring} for intermediate instances with feasible solutions, and normalized dual improvement for hard instances with only dual progress. For easy instances, $\widehat{T}_{\infty}$ estimates the completion time of unfinished runs by fitting a simple surrogate of branch-and-bound tree growth to log snapshots, which places solved and unsolved runs on one time scale. In a retrospective ranking study, it orders $26$ of $34$ configuration pairs correctly, against $17$ for a proxy for the primal--dual integral (Appendix~\ref{app_extrapolation}).

   \item \textit{We show that log-driven semantic tuning improves configuration quality under matched budgets in the easy and intermediate regimes}. We evaluate \name{} on a released benchmark of 215 instances from five benchmark sets, with generated formulations for repeated \textsc{SCIP} tuning and fixed regime labels. In single-instance tuning, \name{} reduces geometric-mean loss by $43.9\%$ and $16.7\%$ in the easy and intermediate regimes, respectively, and achieves lower losses with 10 proposals than the baselines attain with 50. On Case~1, mechanism isolation attributes the gain to the combination of parameter semantics and log diagnostics, while configurations tuned by \name{} on a single instance outperform baseline configurations when transferred to unseen instances of the same family.
\end{itemize}

The remainder of the paper is organized as follows. Section~\ref{sec_related} reviews related work. Section~\ref{sec_problem_formulation} defines the regime-adaptive objective and the log-based completion signal, and Section~\ref{sec_methodology} presents the \name{} framework. Section~\ref{sec_exp} reports the experimental results, and Section~\ref{sec_conclusion} concludes the paper. The appendices provide details on the benchmark, comparison protocol, parameter space, log parser, prompt construction, validation of the completion signal, and additional results. Code, generated benchmark instances, and analysis scripts are available at \coderepo.

\section{Related Work}\label{sec_related}

\paragraph{Automated algorithm configuration.}
Automated solver tuning has traditionally relied on black-box search. Early work demonstrated that automated tuning can outperform expert-designed default settings for MIP solvers \citep{hutter2010automated}. Representative approaches include the SMAC family \citep{hutter2011sequential,lindauer2022smac3}, racing methods such as irace \citep{lopez2016irace}, evolutionary methods such as GGA \citep{ansotegui2009gender}, and surrogate-based optimizers such as RBFOpt \citep{costa2018rbfopt}. These methods treat the solver as an opaque mapping from configuration to outcome and search the configuration space through repeated evaluations. In contrast, LDST uses log-derived diagnostics to identify performance bottlenecks and parameter semantics to guide targeted configuration edits.

\paragraph{Feature-based and instance-specific configuration.}
Classical configurators can also exploit static instance features. SMAC can incorporate instance descriptors into its surrogate model \citep{hutter2011sequential,lindauer2022smac3}; ISAC clusters training instances based on their features and tunes a configuration for each cluster \citep{kadioglu2010isac}; and Hydra-MIP constructs portfolios for heterogeneous MIP families and selects among them for new instances \citep{xu2011hydra}. These approaches use instance features to guide the selection of configurations or portfolios learned from prior instances. In contrast, LDST needs no instance features or configurations learned in advance for a given instance class. Instead, it adapts during tuning using diagnostics extracted from its own solver runs.

\paragraph{Learning-based solver control and configuration.}
Another line of research uses learned models to guide internal solver decisions, including branching \citep{gasse2019exact,gupta2020hybrid} and cut selection \citep{tang2020reinforcement}. Related work on solver configuration uses graph neural networks or Deep Sets to encode problem instances and predict promising parameter settings \citep{hosny2023automatic,song2023instance}. These approaches can exploit problem structure but often require large training sets, closer integration with the solver, or both. Adapting them to new instance distributions may also require retraining. In contrast, LDST operates between solver runs without task-specific training and reads solver behavior from standard logs, avoiding dependence on internal solver states that may be inaccessible without modifying or instrumenting the solver.

\paragraph{Large language models for optimization.}
Recent work has explored LLMs as optimization agents, including OPRO \citep{yang2023large}, heuristic evolution systems such as FunSearch \citep{romera2024mathematical}, and LLM-assisted configuration methods. In hyperparameter optimization, LLMs have been used to propose configurations based on task descriptions and trial histories \citep{zhang2023using,liu2024agenthpo}, while LLAMBO integrates an LLM into Bayesian optimization for warm starting, surrogate modeling, and candidate sampling \citep{liu2024large}. For MIP solvers, LLMs have been used to configure cutting-plane separators \citep{lawless2025llms}, guide online parameter configuration in re-optimization settings \citep{2025online}, and switch among branching rules online based on structured tree-state metrics \citep{2025large}. Execution feedback has also been used to guide LLM-based optimization in other domains, including log-based step-size adaptation for evolution strategies \citep{custode2024investigation}, parallel-program optimization \citep{wei2025interfaces}, and trajectory summaries for hyperparameter tuning \citep{naphade2025expert}. Concurrent work GRIMIP \citep{luo2026grimip} and \name{} both use LLMs for iterative MIP solver configuration without task-specific training. Within a Bayesian optimization framework, GRIMIP draws on the history of evaluated configurations and their scalar performance scores. In contrast, \name{} goes beyond scalar feedback by retaining log-derived diagnostics from successive solver runs and using this process-level history, together with parameter semantics, to guide sparse edits to the best configuration found so far.

\paragraph{Log-based progress extrapolation.}
Recent work has used solver-process signals to estimate search progress and runtime. \citet{hendel2022estimating} study branch-and-bound tree-size estimation using signals such as gap, tree weight, and leaf frequency, together with forecasting and learned estimators. Their results demonstrate that search-completion information can be inferred from solver-process signals, although their strongest estimators rely on richer internal measures and trained predictors implemented within the solver. In contrast, our completion signal uses only the node counts and bounds in the progress table, which MIP solver logs commonly report. By fitting a simple surrogate to the observed trajectory, it estimates completion time using standard log output, without task-specific training, source-code access, or solver modification. In \name{}, this signal places unfinished runs on the same time scale as solved runs, so one objective ranks both on easy instances. The same signal can also support solver monitoring and runtime estimation.

\section{Problem Formulation}
\label{sec_problem_formulation}

In this section, we formalize the tuning setting used throughout the paper. A single scalar loss does not provide equally informative feedback across all instances, because the progress information available at the time limit differs across instances. For instances that the default configuration solves, solved status, final gap, and the primal--dual integral may provide limited information for ranking truncated runs because none directly estimates the remaining time to proof. For instances with feasible solutions but no proof of optimality, both primal and dual bound trajectories are available and can support a trajectory-based measure. When no feasible solution has been found, primal-based measures are unavailable, leaving dual progress as the basis for comparison. We therefore first classify instances according to their outcomes under the default configuration and define a regime-adaptive loss that selects the most informative progress signal available at the time limit. We then introduce a log-based completion signal for easy instances that estimates remaining search effort from the observed branch-and-bound trajectory.

\subsection{Case Split and Regime-Adaptive Loss}
For each instance $i$ and configuration $p$, we run the solver with a fixed time limit $\tau$ (e.g., $600\,\mathrm{s}$), compute a case-dependent scalar loss, and retain the execution log $\Lambda_i(p)$ for diagnostic use. Based on the outcome under the default configuration $p^{\mathrm{def}}$, we assign each instance to one of the three cases summarized in Table~\ref{tab_regime_definition}. Case~1 contains instances solved within $\tau$ under $p^{\mathrm{def}}$ (proved optimal, or proved infeasible for the single infeasible instance). Case~2 contains intermediate instances that are not solved within $\tau$ but do have a primal bound at time $\tau$. Case~3 contains hard instances with a dual bound but no feasible primal solution by $\tau$. These assignments remain fixed throughout tuning, preventing changes in configuration from altering the evaluation objective and ensuring that losses are comparable across methods.

\begin{table}[htbp]
\centering
\caption{Case definitions and tuning signals.}
\label{tab_regime_definition}
\footnotesize
\setlength{\tabcolsep}{3pt}
\renewcommand{\arraystretch}{1.06}
\begin{tabular*}{\linewidth}{@{\extracolsep{\fill}}l l p{0.15\linewidth} p{0.17\linewidth} p{0.29\linewidth}@{}}
\toprule
\textnormal{Case} & \textnormal{Level} & \textnormal{Outcome} & \textnormal{Signal} & \textnormal{Summaries} \\
\midrule
Case~1 & Easy & solved & $\widehat{T}_{\infty}$ & solve rate, time \\
Case~2 & Intermediate & feasible, no proof & PDI & primal, dual, gap \\
Case~3 & Hard & dual progress only & norm.\ dual improv. & final dual, improved-dual coverage \\
\bottomrule
\end{tabular*}
\floatnotes{Cases are defined under the default configuration $p^{\mathrm{def}}$. PDI denotes the primal--dual integral.}
\end{table}

The loss for instance $i$ is defined according to the most informative progress signal available in its assigned case. Let $c_i\in\{1,2,3\}$ denote the case assigned to instance $i$ under the default configuration $p^{\mathrm{def}}$. For Case~1, we use $\widehat{T}_{\infty,i}(p)$, the estimated completion time under configuration $p$, whose construction and fitting procedure are described in Section~\ref{sec_log_completion_signal} and Appendix~\ref{app_extrapolation}. The corresponding normalized loss is $L_i^{(1)}(p)=\widehat{T}_{\infty,i}(p)/\widehat{T}_{\infty,i}(p^{\mathrm{def}})$. Let $z_i^{\mathrm{pr}}(t;p)$ and $z_i^{\mathrm{du}}(t;p)$ denote the primal and dual bounds at time $t$ under configuration $p$. For Case~2, we use the primal--dual integral $\mathrm{PDI}_i(p)$ reported by \textsc{SCIP} at the end of the run \citep{berthold2013measuring}. This quantity integrates the relative gap between $z_i^{\mathrm{pr}}(t;p)$ and $z_i^{\mathrm{du}}(t;p)$ over solving time, with \textsc{SCIP} expressing the gap in percent, capping it at $100$, and using the same value when either bound is unavailable. The corresponding normalized loss is $L_i^{(2)}(p)=\mathrm{PDI}_i(p)/\mathrm{PDI}_i(p^{\mathrm{def}})$. For Case~3, we measure normalized improvement in the final dual bound, defined as $L_i^{(3)}(p)=1+s_i\bigl(z_i^{\mathrm{du}}(\tau;p)-z_i^{\mathrm{du}}(\tau;p^{\mathrm{def}})\bigr)/\bigl(|z_i^{\mathrm{du}}(\tau;p^{\mathrm{def}})|+\varepsilon\bigr)$ with $\varepsilon=10^{-6}$, where $s_i=-1$ for minimization and $s_i=+1$ for maximization. The overall instance-level loss is $L_i(p)=L_i^{(c_i)}(p)$, with every normalized loss floored at $10^{-12}$ to ensure that geometric means remain well defined. This regime-adaptive design selects the strongest reliable progress signal available at the time limit while preserving the interpretability of standard solver summaries.

\subsection{Log-Based Completion Signal}
\label{sec_log_completion_signal}
On Case~1 instances, solve time is the natural objective, but a run that reaches the time limit $\tau$ has no observed solve time. Recording $\tau$ for every such run would tie configurations that differ substantially in their remaining search effort, and final gaps and primal--dual integrals cannot be compared with the solve times of finished runs. We therefore use the completion-time estimator $\widehat{T}_{\infty,i}(p)$, which places finished and unfinished runs on one time scale. If a run is solved within $\tau$, then $\widehat{T}_{\infty,i}(p)=T_i(p)$, where $T_i(p)$ denotes the observed solve time. Otherwise, we truncate the parsed SCIP log at $\tau'=\min\{\tau,t_{\max}\}$, where $t_{\max}$ is the last logged time, and form snapshots $(t_j,n_j,u_j,g_j)$ from its progress-table rows. Here $t_j$ is the elapsed time, $n_j$ the number of processed nodes, $u_j$ the number of open nodes, and $g_j=\bigl|z_{\mathrm{final}}^{\mathrm{pr}}-z^{\mathrm{du}}(t_j)\bigr|$ the absolute difference between the dual bound at $t_j$ and the primal bound $z_{\mathrm{final}}^{\mathrm{pr}}$ reported at the end of the observed run, which for tuning runs is the primal bound available at $\tau$. Following the single-variable branching model of \citet{le2017abstract}, as further analyzed by \citet{anderson2021further}, each snapshot yields the terminal-tree-size surrogate $\widehat{N}_j(x)=n_j+u_jx^{g_j}$. The estimator fits a common terminal tree size $b^\star$ on the log scale, estimates the recent time per processed node $\theta$ from the last five sampled rows, and predicts $\widehat{T}_{\infty,i}(p)=\tau'+\theta\max\{b^\star-n_{\mathrm{last}},0\}$, where $n_{\mathrm{last}}$ is the last sampled processed-node count. The signal is designed for Case~1, where runs are close to proof, and Cases~2 and~3 use the bound-based signals defined above. Appendix~\ref{app_extrapolation} gives the row filtering, the fitting problem, and the safeguards against unstable fits, together with a retrospective validation.

\section{Log-Driven Semantic Tuning}
\label{sec_methodology}

\name{} combines log parsing, sparse edit generation around the current best configuration, and advice validation. Let $\mathcal I_{\mathrm T}$ denote the tuning set and $L(p;\mathcal I_{\mathrm T})$ the aggregate tuning objective, defined as the geometric mean of the instance-level losses in Section~\ref{sec_problem_formulation} over $\mathcal I_{\mathrm T}$. At trial $k$, evaluating configuration $p^{(k)}$ produces aggregate loss $L^{(k)}:=L(p^{(k)};\mathcal I_{\mathrm T})$ and log stream $\Lambda^{(k)}$, while $p^\star$ denotes the best configuration found so far. We first describe the overall tuning loop, followed by the diagnostic and prompting design.

\subsection{The Tuning Loop}
\begin{figure*}[htbp]
  \centering
  \caption{Overview of the \name{} tuning loop.}
  \label{fig:ldst_overview}
  \includegraphics[width=0.9\textwidth]{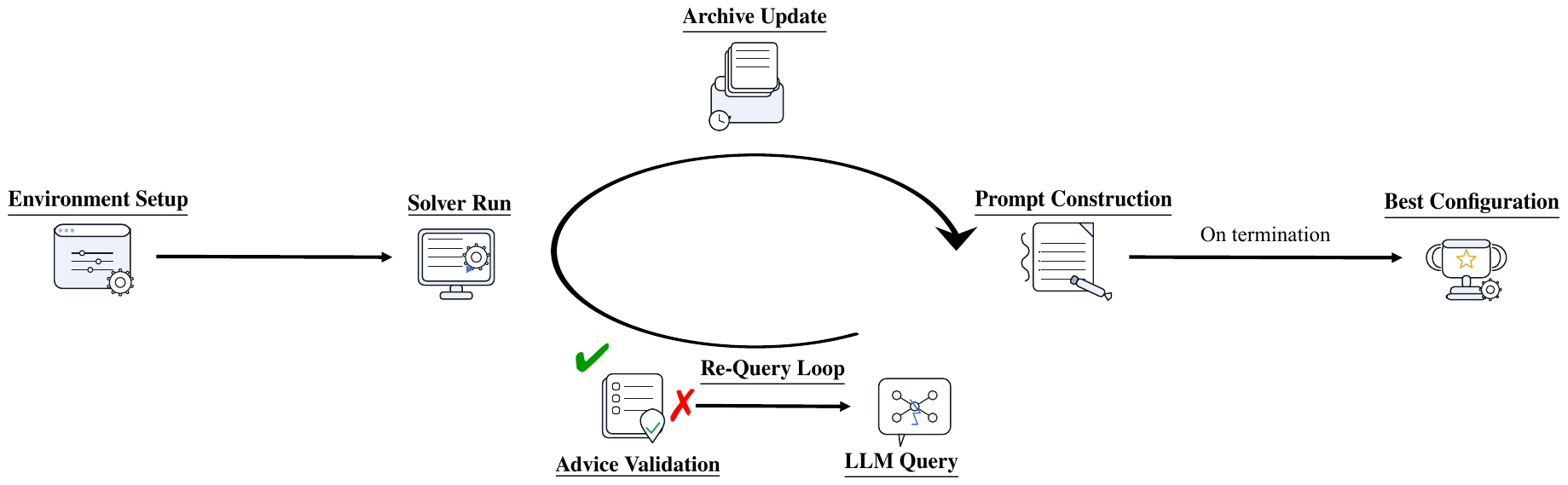}
\floatnotes{Each solver run is parsed and added to the archive. The next prompt is constructed from the current best configuration, bounded memory, and edit constraints. The LLM's proposed edits are validated, with re-queries following failed validation, before the next solver run. Upon termination, the best configuration found is returned.}
\end{figure*}

Figure~\ref{fig:ldst_overview} illustrates the \name{} tuning loop. Starting from the default configuration, each trial runs the solver on $\mathcal I_{\mathrm T}$, extracts diagnostics from the resulting logs, and records the results in an archive. The best configuration $p^\star$ is updated to the archived configuration with the lowest loss, with exact ties broken in favor of the most recent trial. If tuning budget remains, the next prompt is constructed from the current best configuration, bounded memory, and hard edit constraints. The LLM then proposes sparse edits relative to $p^\star$. Advice validation then checks each proposal before execution, and invalid proposals trigger feedback and re-queries until a valid edit is obtained or the re-query limit is reached. The validated configuration is then evaluated in the next trial. By centering the search on the best configuration found so far, \name{} maintains a local and stable search while allowing the LLM to adapt to recent solver behavior.

\subsection{Diagnostic State and Prompting}

\paragraph{Diagnostic summary.}
At trial $k$, the raw log stream $\Lambda^{(k)}$ is summarized in a diagnostic record $r^{(k)}=\bigl(L^{(k)},p^{(k)},\mathbf d^{(k)}\bigr)$, where $\mathbf d^{(k)}=\operatorname{Parse}\!\left(\Lambda^{(k)}\right)$. The diagnostic summary $\mathbf d^{(k)}$ captures search scale, primal-solution and bound progress, and activity in presolving, root relaxation, strong branching, separation, and major heuristics. Applying the same parser and fixed-key schema to every run ensures comparability across trials without instance-specific prompt engineering. This provides \name{} with a consistent view of solver activity and the bottlenecks that dominate each run. Appendix~\ref{app_scip_diagnostics} lists the observable fields used by the parser.

\paragraph{Memory and validation.}
The prompt combines task constraints, a parameter whitelist $\mathcal W$, the current best configuration, and a bounded memory comprising the best historical records and recent trials. Each non-best configuration $p$ is represented relative to $p^\star$ by the edit set $\Delta(p\,;\,p^\star):=\{(j,p_j): j\in\mathcal W,\; p_j\neq p^\star_j\}$. This representation includes only the whitelist parameters whose values differ from the current best configuration, rather than the full parameter vector. Only the entries in $\Delta(p\,;\,p^\star)$ are changed; all other parameter values are inherited from the current best configuration. For example, if a previous configuration differs from the current best configuration only in \texttt{relpscost/minreliable} and \texttt{sbiterquot}, then $\Delta(p\,;\,p^\star)$ contains only these two entries. Before execution, each proposed edit must pass schema, domain, and novelty checks and load successfully as a \textsc{SCIP} settings file. The novelty check uses a deduplication list of previously evaluated configurations. If validation fails, the resulting feedback is included in the next query. When a family-level formulation text is available, the prompt also includes it (Appendix~\ref{app_extra_baselines}). Appendix~\ref{app_prompt_details} provides representative prompts and implementation details.

\section{Experiments}\label{sec_exp}

The experiments test whether log-driven semantic tuning improves single-instance configuration under matched solver-run budgets, identify which information channels contribute to its gains, and measure how well the tuned configurations transfer to unseen instances within the same benchmark family. Exp.~A compares \name{} with \textsc{SMAC3}, \textsc{RBFOpt}, and random search (\textsc{Random}) across the three regimes. Exp.~B isolates the contributions of scalar history, parsed diagnostics, and parameter semantics on Case~1 instances. Exp.~C evaluates within-family transfer from single-instance tuning and multi-instance training splits. All methods share the same solver setup, time limit, editable parameter whitelist, seed convention, and solver-run budgets within each comparison. Appendix~\ref{app_extrapolation} validates $\widehat{T}_{\infty}$, while the remaining appendices provide implementation details, parser fields, and representative prompts.

\subsection{Settings}

\paragraph{Solver setup and configuration space.}
We use \textsc{SCIP} 9.2.4 \citep{bolusani2024scip} in single-thread mode with a per-run time limit of $\tau=600\,\mathrm{s}$. For reproducibility, we fix \texttt{randomization/randomseedshift} to 42 in all solver runs and set the optimizer seeds of \textsc{SMAC3}, \textsc{RBFOpt}, and \textsc{Random} to 42. Each method runs once per instance under this convention, and the sign tests and bootstrap bands below treat instances as the unit of replication. All methods use the same curated whitelist $\mathcal{W}$ of $45$ editable \textsc{SCIP} parameters, detailed in Appendix~\ref{app_whitelist}. We compare \name{} with three black-box configurators, \textsc{SMAC3}, \textsc{RBFOpt}, and \textsc{Random}. All LLM-based methods use GPT-5.2 with medium reasoning effort, and Appendix~\ref{app_hparams} provides the remaining LLM settings and default prompt hyperparameters.

\paragraph{Instances and regime assignment.}
The benchmark comprises 215 instances from five benchmark sets, covering capacitated vehicle routing (CVRP), generalized assignment (GAP), job shop scheduling (JSP), uncapacitated facility location (UFL), and general MIPLIB~2017 instances. Appendix~\ref{app_instance} provides the instance lists and construction details. Following Section~\ref{sec_problem_formulation}, each instance is assigned to Case~1, 2, or~3 based on its outcome under the default configuration $p^{\mathrm{def}}$ at the time limit $\tau$. These assignments remain fixed across all experiments, regardless of the configurations selected during tuning. Table~\ref{tab:exp_regime_stats} summarizes the resulting distribution. Appendix~\ref{app_artifact} documents the code and data release, licensing information, and hardware and software environment used in the experiments.

\begin{table*}[htbp]
\centering
\caption{Case counts by benchmark set with source references.}
\label{tab:exp_regime_stats}
\footnotesize
\setlength{\tabcolsep}{3pt}
\renewcommand{\arraystretch}{1.03}
\begin{tabular*}{\textwidth}{@{\extracolsep{\fill}}llcccc@{\hspace{8pt}}llcccc@{}}
\toprule
\textnormal{Set} & \textnormal{Source} & \textnormal{\#} & \textnormal{Case 1} & \textnormal{Case 2} & \textnormal{Case 3} &
\textnormal{Set} & \textnormal{Source} & \textnormal{\#} & \textnormal{Case 1} & \textnormal{Case 2} & \textnormal{Case 3} \\
\midrule
CVRP   & \cite{lima2014cvrplib}   &  90 & 11 & 39 & 40 & UFL    & \cite{hoefer2018ufllib}   &  30 & 30 &  0 &  0 \\
GAP    & \cite{beasley1990library}&  12 &  7 &  5 &  0 & MIPLIB & \cite{gleixner2021miplib} &  59 & 21 & 30 &  8 \\
JSP    & \cite{beasley1990library}&  24 &  7 & 17 &  0 & \multicolumn{6}{c@{}}{} \\
\midrule
\multicolumn{2}{l}{\textnormal{Total}}          & 215 & 76 & 91 & 48 & \multicolumn{6}{c@{}}{} \\
\bottomrule
\end{tabular*}
\end{table*}

\subsection{Exp. A: Single-Instance Tuning}

\paragraph{A1: Main results.}
We tune each instance independently for $K=50$ proposal trials following an initial default run, yielding $51$ solver evaluations per method-instance pair. The regime-adaptive loss defined in Section~\ref{sec_problem_formulation}, including the Case~1 completion signal $\widehat{T}_{\infty}$, is fixed before comparison and used identically by \name{}, \textsc{SMAC3}, \textsc{RBFOpt}, and \textsc{Random}. For each method, the reported loss on instance $i$ is the best-so-far value $b_i(50)=\min_{0\le t\le 50}L_i(p_i^{(t)})$, where trial~$0$ corresponds to the default configuration. Thus, the reported loss never exceeds~$1$. Table~\ref{tab:exp_main_results} shows that \name{} achieves the lowest geometric-mean loss in Cases~1 and~2, outperforming the best scalar-feedback baseline by $43.9\%$ and $16.7\%$, respectively. These improvements are also reflected in standard solver performance measures. In Case~1, the shifted geometric-mean solve time decreases from $17.1$--$17.9\,\mathrm{s}$ for the baselines to $10.5\,\mathrm{s}$ for \name{}. In Case~2, the final gap decreases from $10.0$--$10.8\%$ to $4.9\%$, accompanied by better primal and dual bounds. Paired statistical tests further support these findings. For each pairwise comparison, we apply an exact two-sided sign test after excluding ties. Let $n_+$ and $n_-$ denote the numbers of instances on which \name{} attains a lower and a higher loss, respectively. The resulting $p$-value equals $\min\{1,\,2\sum_{j=0}^{\min\{n_+,n_-\}}\binom{n_++n_-}{j}2^{-(n_++n_-)}\}$. In Case~1, \name{} outperforms \textsc{SMAC3}, \textsc{RBFOpt}, and \textsc{Random} on $67$, $65$, and $67$ of the $76$ instances, respectively, with all three comparisons significant at $p\le 1.9\times 10^{-10}$. In Case~2, \name{} wins on $71$, $71$, and $68$ of the $91$ instances, with corresponding $p$-values of $7.2\times 10^{-8}$, $7.2\times 10^{-8}$, and $2.5\times 10^{-6}$.

In Case~3, solver runs provide only dual-bound progress within the time limit, and all four methods achieve aggregate losses within $1.2\%$ of the default. \textsc{RBFOpt} attains the lowest aggregate loss ($0.988$), mainly through its final proposal on \texttt{miplib\_irish-electricity}, which improves the dual bound by $30\%$ of its default value (Appendix~\ref{app_sensitivity}). \name{}, \textsc{SMAC3}, \textsc{RBFOpt}, and \textsc{Random} improve the dual bound on $43$, $43$, $41$, and $42$ of the $48$ instances, respectively.

\begin{table*}[htbp]
\centering
\caption{Single-instance tuning results after 50 proposals.}
\label{tab:exp_main_results}
\footnotesize
\setlength{\tabcolsep}{3pt}
\renewcommand{\arraystretch}{1.08}
\begin{tabular*}{\textwidth}{@{\extracolsep{\fill}}l ccc c ccc@{}}
\toprule
& \multicolumn{3}{c}{Primary loss $L$ $\downarrow$}
& \multicolumn{1}{c}{Case~1}
& \multicolumn{3}{c}{Case~2} \\
\cmidrule(lr){2-4}\cmidrule(lr){5-5}\cmidrule(l){6-8}
\textnormal{Method}
& Case~1 & Case~2 & Case~3
& Time (s) $\downarrow$
& Norm. primal $\downarrow$ & Norm. dual $\uparrow$ & Gap (\%) $\downarrow$ \\
\midrule
\textsc{\name{}} & \textbf{0.274} & \textbf{0.557} & 0.993 & \textbf{10.5} & \textbf{1.010} & \textbf{0.990} & \textbf{4.9} \\
\textsc{SMAC3}   & 0.499 & 0.668 & 0.993 & 17.3 & 1.034 & 0.975 & 10.3 \\
\textsc{RBFOpt}  & 0.490 & 0.675 & \textbf{0.988} & 17.1 & 1.039 & 0.975 & 10.8 \\
\textsc{Random}  & 0.519 & 0.669 & 0.993 & 17.9 & 1.030 & 0.977 & 10.0 \\
\bottomrule
\end{tabular*}
\floatnotes{For each method, the selected configuration is the best among the default (trial~0) and the 50 proposals. Primary losses are geometric means over all 76, 91, and 48 instances in Cases~1, 2, and~3, respectively. ``Time'' is the shifted geometric mean of solve time (shift $1\,\mathrm{s}$) over the 75 Case~1 instances that every method solves to optimality; the remaining instance, \texttt{miplib\_neos859080}, is infeasible. For Case~2, let $\tilde z^{\mathrm{pr}}_{i,m}=-s_i z^{\mathrm{pr}}_{i,m}$ and $\tilde z^{\mathrm{du}}_{i,m}=s_i z^{\mathrm{du}}_{i,m}$ denote the sense-aligned final primal and dual bounds of method $m$ on instance $i$, and let $b_i^\star$ denote the best aligned bound among the compared methods on the 90 instances with finite bounds under all methods. ``Norm. primal'' is the geometric mean of $1+(\tilde z^{\mathrm{pr}}_{i,m}-b_i^\star)/\max\{|\tilde z^{\mathrm{pr}}_{i,m}|,|b_i^\star|,\varepsilon\}$ (lower is better). ``Norm. dual'' applies the reciprocal score to the aligned dual bound (higher is better). ``Gap'' is the shifted geometric mean (shift $0.01$) of the final relative gap in percent on the same 90 instances.}
\end{table*}

\begin{figure*}[htbp]
  \centering
  \caption{Geometric mean of the best-so-far loss in Exp.~A.}
  \label{fig:exp_a1_curves}
  \includegraphics[width=0.98\linewidth]{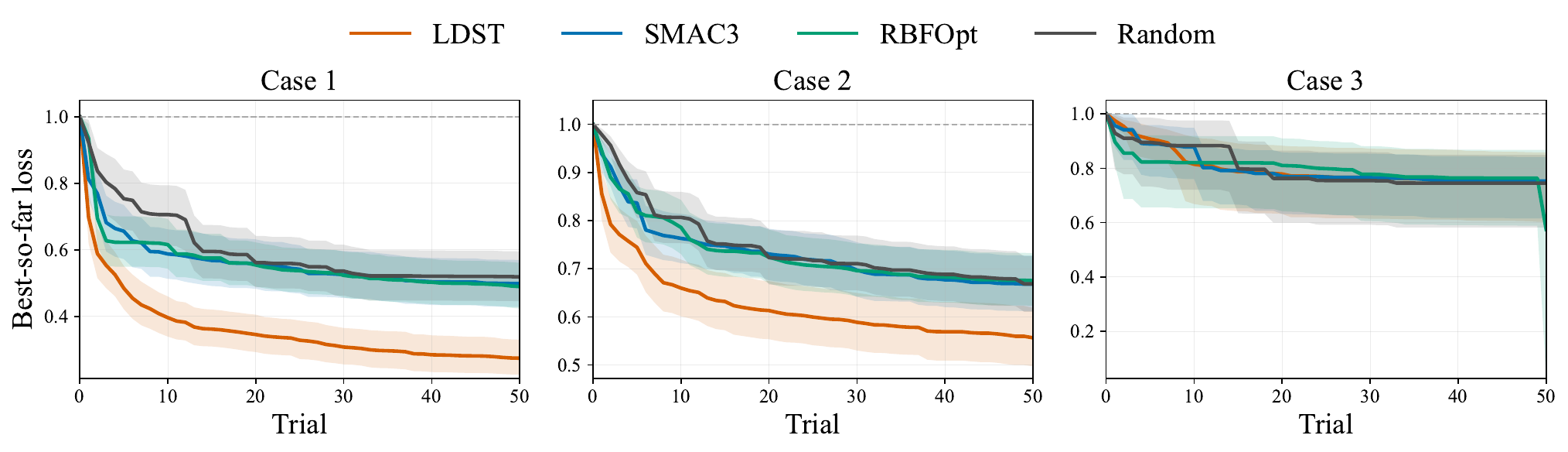}
\floatnotes{For each regime $c$ and trial $t$, the shaded band is constructed from $2000$ bootstrap samples, each containing $|\mathcal I_c|$ instances drawn with replacement from $\mathcal I_c$. For each sample, the aggregate best-so-far loss at trial $t$ is recomputed using the same formula as the main curve. The lower and upper edges are the 2.5th and 97.5th percentiles of the $2000$ bootstrap values at that trial. All curves start from the default configuration at trial~0 (loss~1). The Case~3 panel applies the display-only transformation $\widetilde{L}^{(3)}=1+36(L^{(3)}-1)$.}
\end{figure*}

\paragraph{A2: Convergence behavior.}
Figure~\ref{fig:exp_a1_curves} shows that \name{} establishes an early advantage and maintains it throughout tuning in Cases~1 and~2. In Case~1, \name{} achieves the lowest geometric-mean best-so-far loss of $0.396$ by trial~10, compared with $0.588$--$0.706$ for the baselines. By trial~50, its loss decreases to $0.274$, compared with $0.490$ for \textsc{RBFOpt}, $0.499$ for \textsc{SMAC3}, and $0.519$ for \textsc{Random}. A similar pattern emerges in Case~2. By trial~10, \name{} reaches a loss of $0.660$, while the baselines remain at $0.763$--$0.807$. By trial~50, its loss decreases to $0.557$, compared with $0.668$--$0.675$ for the three baselines. Thus, in both regimes, \name{} achieves a lower loss with 10 proposals than any baseline attains with 50. Exp.~B identifies which information channels produce this early advantage.

For Case~3, Figure~\ref{fig:exp_a1_curves} applies the display-only transformation $\widetilde{L}^{(3)}=1+36(L^{(3)}-1)$ to make the small differences between methods visible. The late decrease in the \textsc{RBFOpt} curve reflects its trial-50 improvement on \texttt{miplib\_irish-electricity}, as discussed in A1.

\paragraph{A3: Edit sparsity and LLM cost.}
\name{} changes few parameters per proposal. Across its $10{,}750$ proposals in Exp.~A, the evaluated configuration differs from the best configuration at the time of the proposal in a median of $7$ of the $45$ whitelisted parameters, with medians of $5$, $7$, and $11$ in Cases~1, 2, and~3. Overall, $82\%$ of the proposals change at most $10$ parameters, and $1.2\%$ reach the edit cap of $20$ (Table~\ref{tab:edit_size} in Appendix~\ref{app_sensitivity}). The LLM side is also inexpensive. The $10{,}964$ LLM calls behind these proposals, including re-queries, cost about USD~$320$ at uncached list prices, or about USD~$0.03$ per proposal (Table~\ref{tab:artifact_env} in Appendix~\ref{app_artifact}).

\subsection{Exp. B: Mechanism Isolation}\label{subsec:exp_log_ablation}

Exp.~B examines the information channels underlying the gains of \name{}. We compare \name{} with \textsc{LLM-StaticContext}, \textsc{LLM-OutcomeOnly}, and \textsc{LLM-IndexOnly} under the single-instance Case~1 protocol. Each of the 76 Case~1 instances is tuned independently for $K=30$ trials, and results are aggregated over all instances. The \name{} results use the first 30 trials of its Exp.~A runs. All methods follow the shared execution contract in Appendix~\ref{app_extra_baselines}, with the controls differing only in the information available to the decision policy. Specifically, \textsc{LLM-StaticContext} removes the cross-trial archive, including scalar history, along with the diagnostic guide and deduplication list, while retaining the parsed diagnostics of the current best run. \textsc{LLM-OutcomeOnly} retains scalar history but removes parsed diagnostics and the diagnostic guide. \textsc{LLM-IndexOnly} replaces parameter names and defaults with indices and domain hints. Table~\ref{tab:prompt_ablation_blocks} details these changes to the prompts.

Table~\ref{tab:exp_log_ablation_main} shows that \name{} performs best across all three metrics, achieving a Best loss of $0.308$, early aggregate loss $L_{\le 5}=0.485$, and threshold-hit rate $H_{0.8}=90.8\%$. All three metrics give the same ranking, with \name{} followed by \textsc{LLM-StaticContext}, \textsc{LLM-OutcomeOnly}, and \textsc{LLM-IndexOnly}. Removing parameter semantics causes the largest performance degradation. Despite retaining scalar history and parsed diagnostics, \textsc{LLM-IndexOnly} achieves a Best loss of $0.523$, with $L_{\le 5}$ at $0.707$ and $H_{0.8}$ at $76.3\%$. This suggests that index-level domain hints alone are insufficient to translate solver feedback into effective parameter edits. Removing scalar history and archive memory has a smaller but noticeable effect. Although \textsc{LLM-StaticContext} retains the current best run's diagnostics, it remains weaker than \name{}, with a Best loss of $0.342$ versus $0.308$, $0.509$ versus $0.485$ in $L_{\le 5}$, and $86.8\%$ versus $90.8\%$ in $H_{0.8}$. Removing parsed diagnostics while retaining scalar history also reduces performance. \textsc{LLM-OutcomeOnly} achieves a Best loss of $0.358$, $L_{\le 5}=0.520$, and $H_{0.8}=84.2\%$, a smaller degradation than that caused by removing parameter semantics. All three controls therefore underperform the full method. Parameter semantics contribute the most, while parsed diagnostics provide a larger benefit than cross-trial memory. These results indicate that the advantage of \name{} comes from combining knowledge of parameter semantics with diagnostic evidence from solver logs.

Parameter semantics also help maintain a productive search. Without parameter names and defaults, \textsc{LLM-IndexOnly} exhausts the re-query limit before trial~30 on $6$ of the $76$ instances, five times by proposing a previously evaluated configuration and once by returning no edits. \textsc{LLM-StaticContext}, which lacks the archive and the deduplication list, terminates early on one instance, whereas \name{} completes all trials on every instance. For early-terminated runs, the last best-so-far loss is carried forward to $b_i(K)$.

\begin{table}[htbp]
\centering
\caption{Mechanism isolation in Exp.~B on the 76 Case~1 instances.}
\label{tab:exp_log_ablation_main}
\scriptsize
\setlength{\tabcolsep}{2pt}
\renewcommand{\arraystretch}{1.05}
\begin{tabular*}{\columnwidth}{@{\extracolsep{\fill}}lcccccc@{}}
\toprule
\textnormal{Method}
& \textnormal{\shortstack{Scalar\\history}}
& \textnormal{\shortstack{Parsed\\diagnostics}}
& \textnormal{\shortstack{Parameter\\semantics}}
& \textnormal{Best loss} $\downarrow$
& $L_{\le 5}$ $\downarrow$
& $H_{0.8}$ (\%) $\uparrow$ \\
\midrule
\textsc{\name{}} & \checkmark & \checkmark & \checkmark & \textbf{0.308} & \textbf{0.485} & \textbf{90.8} \\
\textsc{LLM-StaticContext} & \xmark & \checkmark & \checkmark & 0.342 & 0.509 & 86.8 \\
\textsc{LLM-OutcomeOnly} & \checkmark & \xmark & \checkmark & 0.358 & 0.520 & 84.2 \\
\textsc{LLM-IndexOnly} & \checkmark & \checkmark & \xmark & 0.523 & 0.707 & 76.3 \\
\bottomrule
\end{tabular*}
\floatnotes{The access columns indicate whether scalar history, parsed diagnostics, and parameter semantics are retained relative to \name{}. Index-level domain hints in \textsc{LLM-IndexOnly} are not counted as parameter semantics. Let $b_i(t):=\min_{0\le s\le t}L_i(p_i^{(s)})$, and let $\mathrm{GM}$ denote the geometric mean. Best loss is $\mathrm{GM}_{i\in\mathcal I}b_i(K)$, early aggregate loss is $L_{\le 5}:=\mathrm{GM}_{i\in\mathcal I}b_i(5)$, and the threshold-hit rate is $H_{0.8}:=100|\mathcal I|^{-1}\sum_{i\in\mathcal I}\mathbf{1}\{b_i(K)\le0.8\}$. For early-terminated runs, the last best-so-far loss is carried forward.}
\end{table}

\subsection{Exp. C: Within-Family Transfer from Single and Multiple Training Instances}\label{subsec:exp_generalization}

\paragraph{C1: Transfer from a single anchor.}
C1 examines whether a configuration tuned on a single instance transfers to other instances within the same benchmark family. We use the 55 Case~1 instances from CVRP, GAP, JSP, and UFL as anchors. For each anchor and method, we select the configuration obtained in Exp.~A and evaluate it, without further tuning, on the remaining Case~1 instances from the same family. \name{} improves upon the default configuration on all $55$ anchors, whereas \textsc{SMAC3}, \textsc{RBFOpt}, and \textsc{Random} select the default on $4$, $3$, and $4$ anchors, respectively. Table~\ref{tab:exp_one_instance_transfer} shows that \name{} achieves the lowest anchor and test losses overall and the lowest test loss in every family. Its aggregate test loss is $0.499$, compared with $0.798$--$0.838$ for the baselines. Paired sign tests across the $55$ anchors further support this advantage, with \name{} outperforming \textsc{SMAC3}, \textsc{RBFOpt}, and \textsc{Random} on $48$, $49$, and $54$ anchors, respectively (all $p\le 1.4\times 10^{-8}$).

\begin{table*}[htbp]
\centering
\caption{C1 one-instance transfer on the 55 Case~1 anchors.}
\label{tab:exp_one_instance_transfer}
\footnotesize
\setlength{\tabcolsep}{3pt}
\renewcommand{\arraystretch}{1.05}
\begin{tabular*}{\textwidth}{@{\extracolsep{\fill}}lccccccc@{}}
\toprule
& \multicolumn{3}{c}{\textnormal{Overall}} & \multicolumn{4}{c}{\textnormal{Family-level Case~1 test}} \\
\cmidrule(lr){2-4}\cmidrule(lr){5-8}
\textnormal{Method} & \textnormal{Case~1 train} & \textnormal{Case~1 test} & \textnormal{Test$-$train} & \textnormal{CVRP} & \textnormal{GAP} & \textnormal{JSP} & \textnormal{UFL} \\
\midrule
\textsc{\name{}} & \textbf{0.238} & \textbf{0.499} & 0.261 & \textbf{1.190} & \textbf{0.734} & \textbf{0.262} & \textbf{0.385} \\
\textsc{SMAC3} & 0.429 & 0.799 & 0.370 & 2.067 & 0.966 & 0.657 & 0.565 \\
\textsc{RBFOpt} & 0.427 & 0.798 & 0.370 & 1.978 & 1.078 & 0.716 & 0.547 \\
\textsc{Random} & 0.449 & 0.838 & 0.390 & 1.866 & 1.313 & 0.635 & 0.601 \\
\bottomrule
\end{tabular*}
\floatnotes{Each anchor configuration is selected as the best among the default configuration and the 50 proposals in Exp.~A for that anchor. Entries are geometric means over anchor splits. Test$-$train denotes $L_{\mathrm{test}}-L_{\mathrm{train}}$. The last four columns report Case~1 test loss separately for CVRP, GAP, JSP, and UFL.}
\end{table*}

\paragraph{C2: Multi-instance tuning versus independent selection.}
C2 compares two ways to use multiple training instances from the same family, independent selection and joint tuning. For each eligible benchmark set $s$ and training size $m\in\{2,3,4,5,6\}$, we construct five deterministic Case~1 train/test splits, provided that $s$ contains at least $m+4$ Case~1 instances, leaving at least four instances for testing. Independent selection uses the configurations obtained from the individual Exp.~A runs and selects the one with the lowest geometric-mean loss across the training subset. Joint tuning instead runs \name{} for 50 proposals on the full training subset. The selected configurations are then evaluated on the Case~1 test instances and, when available, on Case~2 and Case~3 instances from the same family.

Table~\ref{tab:exp_gen_summary} compares the two strategies. With two or three training instances, joint tuning achieves lower Case~1 training and test losses and lower Case~2 and Case~3 losses. Its stronger fit to the training subset increases the train--test gap from $0.150$ to $0.226$ at $m=2$ and from $0.099$ to $0.143$ at $m=3$, even as the test loss decreases. At $m=4$, the two strategies achieve nearly identical Case~1 test losses ($0.560$ for independent selection versus $0.561$ for joint tuning), although joint tuning still yields lower training loss and lower losses in Case~2 and Case~3. For larger training subsets, the relative performance depends on the benchmark family. Joint tuning outperforms independent selection on CVRP at every training size and on GAP at both eligible sizes, whereas independent selection performs better on UFL. Consequently, the pooled averages depend on the composition of eligible benchmark families at each training size (Appendix~\ref{app_sensitivity}).

\begin{table}[H]
\centering
\caption{Within-family transfer under independent and joint tuning in C2.}
\label{tab:exp_gen_summary}
\footnotesize
\setlength{\tabcolsep}{1pt}
\renewcommand{\arraystretch}{1.02}
\begin{tabular*}{\textwidth}{@{\extracolsep{\fill}}c cc cc cc cc cc@{}}
\toprule
\textnormal{$m$}
& \multicolumn{2}{c}{\textnormal{Case~1 train}}
& \multicolumn{2}{c}{\textnormal{Case~1 test}}
& \multicolumn{2}{c}{\textnormal{Test$-$train $\downarrow$}}
& \multicolumn{2}{c}{\textnormal{Case~2}}
& \multicolumn{2}{c}{\textnormal{Case~3}} \\
\cmidrule(lr){2-3}\cmidrule(lr){4-5}\cmidrule(lr){6-7}\cmidrule(lr){8-9}\cmidrule(lr){10-11}
& \textnormal{Indep.} & \textnormal{Joint}
& \textnormal{Indep.} & \textnormal{Joint}
& \textnormal{Indep.} & \textnormal{Joint}
& \textnormal{Indep.} & \textnormal{Joint}
& \textnormal{Indep.} & \textnormal{Joint} \\
\midrule
$2$ & 0.390 & \textbf{0.282} & 0.540 & \textbf{0.508} & \textbf{0.150} & 0.226 & 0.919 & \textbf{0.849} & 1.010 & \textbf{1.002} \\
$3$ & 0.396 & \textbf{0.325} & 0.495 & \textbf{0.468} & \textbf{0.099} & 0.143 & 0.864 & \textbf{0.828} & 1.018 & \textbf{1.011} \\
$4$ & 0.421 & \textbf{0.395} & \textbf{0.560} & 0.561 & \textbf{0.139} & 0.166 & 0.715 & \textbf{0.689} & 1.025 & \textbf{1.007} \\
$5$ & \textbf{0.389} & 0.435 & \textbf{0.543} & 0.625 & \textbf{0.154} & 0.190 & \textbf{0.687} & 0.775 & 1.030 & \textbf{1.008} \\
$6$ & 0.467 & \textbf{0.376} & \textbf{0.474} & 0.508 & \textbf{0.007} & 0.132 & \textbf{0.705} & 0.720 & 1.009 & \textbf{1.008} \\
\bottomrule
\end{tabular*}
\floatnotes{Each metric is reported for both Independent and Joint tuning. Test$-$train denotes $L_{\mathrm{test}}-L_{\mathrm{train}}$. At $m=2,3$, Case~1 includes all four benchmark families, Case~2 includes CVRP, GAP, and JSP, and Case~3 includes only CVRP. At $m=4,5,6$, Case~1 includes CVRP and UFL, while Cases~2 and~3 include only CVRP.}
\end{table}

\section{Conclusion}\label{sec_conclusion}

We introduced \name{}, which uses solver-log diagnostics and an LLM to propose validated sparse edits to the current best configuration. The framework combines memory centered on that configuration, a regime-adaptive loss, and a Case~1 completion-time signal based on branch-and-bound tree growth. Across 215 MIP instances, \name{} achieves final geometric-mean losses of $0.274$ and $0.557$ in the easy and intermediate regimes, respectively, $43.9\%$ and $16.7\%$ below the best scalar-feedback baselines. In both regimes, \name{} achieves with 10 proposals a lower loss than any baseline attains with 50. On Case~1, the mechanism study attributes the gain to parameter semantics combined with log diagnostics, and configurations that \name{} tunes on a single instance transfer better than baseline configurations to unseen instances of the same family. Semantic reasoning over these diagnostics can make MIP configuration substantially more sample efficient. The reliance on ordinary logs and parameter names suggests extending \name{} to other solvers, larger configuration spaces, and other language models.

\bibliographystyle{informs2014}
\bibliography{references}

\clearpage
\appendix

\section*{Appendix Roadmap}
Appendix~\ref{app_instance} describes the benchmark sources, MPS-construction protocol, default-based regime assignment, and within-family split design used in the main experiments. Appendix~\ref{app_extra_baselines} specifies the shared execution contract, scalar-feedback baselines, and Exp.~B controls. Appendices~\ref{app_whitelist} and \ref{app_hparams} document the configuration space and default hyperparameters, while Appendices~\ref{app_scip_diagnostics} and \ref{app_prompt_details} describe the fixed log parser and prompt construction used by \name{}. Appendix~\ref{app_extrapolation} validates the Case~1 completion-time estimator $\widehat{T}_{\infty}$. Appendix~\ref{app_sensitivity} provides the edit sizes of Exp.~A, additional Exp.~C results, and an audit of the Case~3 outlier. Appendix~\ref{app_artifact} documents the code and data release, licensing information, hardware and software environment, and token-side compute accounting. Finally, Appendix~\ref{app_instance_lists} provides the complete family-specific formulations and benchmark manifests.

\section{Benchmark Sources, Formulations, and Regime Assignment}\label{app_instance}
This appendix describes the benchmark design underlying the main experiments, including the source collections, MPS-construction protocol, default-based regime assignment, and within-family split design used in Exp.~C. Appendix~\ref{app_instance_lists} provides the complete family-specific formulations and benchmark manifests.

We use five public benchmark sets, CVRPLIB CVRP \citep{lima2014cvrplib}, OR-Library GAP and JSP \citep{beasley1990library}, UflLib UFL \citep{hoefer2018ufllib}, and MIPLIB~2017 \citep{gleixner2021miplib}. The final benchmark comprises $215$ instances, including $90$ CVRP, $12$ GAP, $24$ JSP, $30$ UFL, and $59$ MIPLIB instances. For the four domain-specific families, we export standard mixed-integer linear programming (MILP) formulations to \texttt{.mps} files through an automated pipeline. Each generated file is checked against its native model for variable and constraint counts, integrality patterns, and objective agreement within numerical tolerances. For MIPLIB, we initially select 60 instances from the official \texttt{benchmark-v2} subset using a fixed random seed. We exclude \texttt{k1mushroom} because it does not complete presolving under default SCIP settings within the 600-second time limit, leaving 59 MIPLIB instances for the reported experiments. The code repository provides the final manifests, file-level provenance, and regime labels (Appendix~\ref{app_artifact}).

\subsection{Default-Based Regime Assignment}
Each instance is assigned to a case based on its outcome under the default configuration $p^{\mathrm{def}}$ with a common time limit of $\tau=600\,\mathrm{s}$, following the definitions in Section~\ref{sec_problem_formulation}. Case~1 contains instances solved under $p^{\mathrm{def}}$, including those proved optimal and the single instance \texttt{miplib\_neos859080}, which is proved infeasible. Case~2 contains instances with a feasible primal solution but no proof of optimality, while Case~3 contains the remaining hard instances evaluated using the dual-progress objective. These assignments remain fixed across all subsequent experiments and are not recomputed after tuning. This ensures consistent regime definitions across methods and prevents tuned configurations from changing the evaluation objective.

Table~\ref{tab:benchmark_family_overview} summarizes the source collections, construction procedures, and default-based regime counts for each benchmark family.
\begin{table*}[htbp]
\centering
\caption{Benchmark sources, construction routes, and regime counts.}
\label{tab:benchmark_family_overview}
\footnotesize
\setlength{\tabcolsep}{4pt}
\renewcommand{\arraystretch}{1.06}
\begin{tabular*}{\textwidth}{@{\extracolsep{\fill}}l l l r r r r@{}}
\toprule
\textnormal{Family} & \textnormal{Source collection} & \textnormal{Construction route} & \textnormal{\#} & \textnormal{Case 1} & \textnormal{Case 2} & \textnormal{Case 3} \\
\midrule
CVRP & CVRPLIB & two-commodity flow MILP export & 90 & 11 & 39 & 40 \\
GAP & OR-Library \texttt{gapc/gapd} & assignment MILP export & 12 & 7 & 5 & 0 \\
JSP & OR-Library \texttt{jobshop1} & disjunctive MILP export & 24 & 7 & 17 & 0 \\
UFL & UflLib Uniform & facility-location MILP export & 30 & 30 & 0 & 0 \\
MIPLIB & MIPLIB~2017 \texttt{benchmark-v2} & direct MPS benchmark instances & 59 & 21 & 30 & 8 \\
\bottomrule
\end{tabular*}
\floatnotes{Counts cover the 215-instance benchmark. Regime assignments are made under the default configuration and remain fixed.}
\end{table*}

Table~\ref{tab:case_assignment_manifest} gives the complete per-instance case-assignment manifest used by all experiments. Table~\ref{tab:benchmark_family_overview} summarizes the corresponding counts.
\begin{table}[p]
\centering
\caption{Case assignments for all benchmark instances.}
\label{tab:case_assignment_manifest}
{\footnotesize
\setlength{\tabcolsep}{3pt}
\renewcommand{\arraystretch}{0.92}
\begin{tabular}{@{}L{0.076\textwidth}L{0.335\textwidth}C{0.044\textwidth}@{\hspace{8pt}}L{0.076\textwidth}L{0.335\textwidth}C{0.044\textwidth}@{}}
\toprule
\textnormal{Family} & \textnormal{Instance} & \textnormal{Case} &
\textnormal{Family} & \textnormal{Instance} & \textnormal{Case} \\
\midrule
\textnormal{CVRP} & \texttt{A-n32-k5} & 2 & \textnormal{CVRP} & \texttt{A-n33-k5} & 2 \\
\textnormal{CVRP} & \texttt{A-n33-k6} & 2 & \textnormal{CVRP} & \texttt{A-n34-k5} & 2 \\
\textnormal{CVRP} & \texttt{A-n36-k5} & 2 & \textnormal{CVRP} & \texttt{A-n37-k5} & 2 \\
\textnormal{CVRP} & \texttt{A-n37-k6} & 2 & \textnormal{CVRP} & \texttt{A-n38-k5} & 2 \\
\textnormal{CVRP} & \texttt{A-n39-k5} & 2 & \textnormal{CVRP} & \texttt{A-n39-k6} & 2 \\
\textnormal{CVRP} & \texttt{A-n44-k6} & 2 & \textnormal{CVRP} & \texttt{A-n45-k6} & 3 \\
\textnormal{CVRP} & \texttt{A-n45-k7} & 2 & \textnormal{CVRP} & \texttt{A-n46-k7} & 2 \\
\textnormal{CVRP} & \texttt{A-n48-k7} & 2 & \textnormal{CVRP} & \texttt{A-n53-k7} & 2 \\
\textnormal{CVRP} & \texttt{A-n54-k7} & 3 & \textnormal{CVRP} & \texttt{A-n55-k9} & 3 \\
\textnormal{CVRP} & \texttt{A-n60-k9} & 3 & \textnormal{CVRP} & \texttt{A-n61-k9} & 3 \\
\textnormal{CVRP} & \texttt{A-n62-k8} & 3 & \textnormal{CVRP} & \texttt{A-n63-k10} & 3 \\
\textnormal{CVRP} & \texttt{A-n63-k9} & 3 & \textnormal{CVRP} & \texttt{A-n64-k9} & 3 \\
\textnormal{CVRP} & \texttt{A-n65-k9} & 3 & \textnormal{CVRP} & \texttt{A-n69-k9} & 3 \\
\textnormal{CVRP} & \texttt{A-n80-k10} & 3 & \textnormal{CVRP} & \texttt{B-n31-k5} & 2 \\
\textnormal{CVRP} & \texttt{B-n34-k5} & 2 & \textnormal{CVRP} & \texttt{B-n35-k5} & 2 \\
\textnormal{CVRP} & \texttt{B-n38-k6} & 2 & \textnormal{CVRP} & \texttt{B-n39-k5} & 2 \\
\textnormal{CVRP} & \texttt{B-n41-k6} & 2 & \textnormal{CVRP} & \texttt{B-n43-k6} & 2 \\
\textnormal{CVRP} & \texttt{B-n44-k7} & 2 & \textnormal{CVRP} & \texttt{B-n45-k5} & 2 \\
\textnormal{CVRP} & \texttt{B-n45-k6} & 3 & \textnormal{CVRP} & \texttt{B-n50-k7} & 2 \\
\textnormal{CVRP} & \texttt{B-n50-k8} & 2 & \textnormal{CVRP} & \texttt{B-n51-k7} & 3 \\
\textnormal{CVRP} & \texttt{B-n52-k7} & 2 & \textnormal{CVRP} & \texttt{B-n56-k7} & 2 \\
\textnormal{CVRP} & \texttt{B-n57-k7} & 3 & \textnormal{CVRP} & \texttt{B-n57-k9} & 3 \\
\textnormal{CVRP} & \texttt{B-n63-k10} & 3 & \textnormal{CVRP} & \texttt{B-n64-k9} & 3 \\
\textnormal{CVRP} & \texttt{B-n66-k9} & 3 & \textnormal{CVRP} & \texttt{B-n67-k10} & 3 \\
\textnormal{CVRP} & \texttt{B-n68-k9} & 3 & \textnormal{CVRP} & \texttt{B-n78-k10} & 3 \\
\textnormal{CVRP} & \texttt{E-n101-k14} & 3 & \textnormal{CVRP} & \texttt{E-n101-k8} & 3 \\
\textnormal{CVRP} & \texttt{E-n13-k4} & 1 & \textnormal{CVRP} & \texttt{E-n22-k4} & 1 \\
\textnormal{CVRP} & \texttt{E-n23-k3} & 1 & \textnormal{CVRP} & \texttt{E-n30-k3} & 2 \\
\textnormal{CVRP} & \texttt{E-n31-k7} & 1 & \textnormal{CVRP} & \texttt{E-n33-k4} & 2 \\
\textnormal{CVRP} & \texttt{E-n51-k5} & 2 & \textnormal{CVRP} & \texttt{E-n76-k10} & 3 \\
\textnormal{CVRP} & \texttt{E-n76-k14} & 3 & \textnormal{CVRP} & \texttt{E-n76-k7} & 3 \\
\textnormal{CVRP} & \texttt{E-n76-k8} & 3 & \textnormal{CVRP} & \texttt{F-n135-k7} & 3 \\
\textnormal{CVRP} & \texttt{F-n45-k4} & 2 & \textnormal{CVRP} & \texttt{F-n72-k4} & 2 \\
\textnormal{CVRP} & \texttt{P-n101-k4} & 2 & \textnormal{CVRP} & \texttt{P-n16-k8} & 1 \\
\textnormal{CVRP} & \texttt{P-n19-k2} & 1 & \textnormal{CVRP} & \texttt{P-n20-k2} & 1 \\
\textnormal{CVRP} & \texttt{P-n21-k2} & 1 & \textnormal{CVRP} & \texttt{P-n22-k2} & 1 \\
\textnormal{CVRP} & \texttt{P-n22-k8} & 1 & \textnormal{CVRP} & \texttt{P-n23-k8} & 2 \\
\textnormal{CVRP} & \texttt{P-n40-k5} & 1 & \textnormal{CVRP} & \texttt{P-n45-k5} & 2 \\
\textnormal{CVRP} & \texttt{P-n50-k10} & 3 & \textnormal{CVRP} & \texttt{P-n50-k7} & 2 \\
\textnormal{CVRP} & \texttt{P-n50-k8} & 3 & \textnormal{CVRP} & \texttt{P-n51-k10} & 3 \\
\textnormal{CVRP} & \texttt{P-n55-k10} & 3 & \textnormal{CVRP} & \texttt{P-n55-k15} & 3 \\
\textnormal{CVRP} & \texttt{P-n55-k7} & 2 & \textnormal{CVRP} & \texttt{P-n55-k8} & 2 \\
\textnormal{CVRP} & \texttt{P-n60-k10} & 3 & \textnormal{CVRP} & \texttt{P-n60-k15} & 3 \\
\textnormal{CVRP} & \texttt{P-n65-k10} & 3 & \textnormal{CVRP} & \texttt{P-n70-k10} & 3 \\
\textnormal{CVRP} & \texttt{P-n76-k4} & 3 & \textnormal{CVRP} & \texttt{P-n76-k5} & 3 \\
\textnormal{GAP} & \texttt{gapc-1} & 1 & \textnormal{GAP} & \texttt{gapc-2} & 1 \\
\textnormal{GAP} & \texttt{gapc-3} & 1 & \textnormal{GAP} & \texttt{gapc-4} & 1 \\
\textnormal{GAP} & \texttt{gapc-5} & 1 & \textnormal{GAP} & \texttt{gapc-6} & 1 \\
\textnormal{GAP} & \texttt{gapd-1} & 1 & \textnormal{GAP} & \texttt{gapd-2} & 2 \\
\textnormal{GAP} & \texttt{gapd-3} & 2 & \textnormal{GAP} & \texttt{gapd-4} & 2 \\
\textnormal{GAP} & \texttt{gapd-5} & 2 & \textnormal{GAP} & \texttt{gapd-6} & 2 \\
\textnormal{JSP} & \texttt{abz5} & 2 & \textnormal{JSP} & \texttt{abz6} & 1 \\
\textnormal{JSP} & \texttt{abz7} & 2 & \textnormal{JSP} & \texttt{abz8} & 2 \\
\textnormal{JSP} & \texttt{abz9} & 2 & \textnormal{JSP} & \texttt{ft10} & 2 \\
\bottomrule
\end{tabular}
}
\floatnotes{Each instance is assigned to a case once under the default configuration $p^{\mathrm{def}}$ with time limit $\tau=600\,\mathrm{s}$. These assignments remain fixed across all experiments.}
\end{table}

\begin{table}[p]
\centering
\addtocounter{table}{-1}
\renewcommand{\theHtable}{\arabic{table}.continued}
\caption{Case assignments for all benchmark instances (continued).}
\label{tab:case_assignment_manifest_continued}
{\footnotesize
\setlength{\tabcolsep}{3pt}
\renewcommand{\arraystretch}{0.92}
\begin{tabular}{@{}L{0.076\textwidth}L{0.335\textwidth}C{0.044\textwidth}@{\hspace{8pt}}L{0.076\textwidth}L{0.335\textwidth}C{0.044\textwidth}@{}}
\toprule
\textnormal{Family} & \textnormal{Instance} & \textnormal{Case} &
\textnormal{Family} & \textnormal{Instance} & \textnormal{Case} \\
\midrule
\textnormal{JSP} & \texttt{ft20} & 2 & \textnormal{JSP} & \texttt{la01} & 1 \\
\textnormal{JSP} & \texttt{la02} & 1 & \textnormal{JSP} & \texttt{la03} & 1 \\
\textnormal{JSP} & \texttt{la04} & 1 & \textnormal{JSP} & \texttt{la05} & 1 \\
\textnormal{JSP} & \texttt{la06} & 2 & \textnormal{JSP} & \texttt{la07} & 2 \\
\textnormal{JSP} & \texttt{la08} & 2 & \textnormal{JSP} & \texttt{la09} & 2 \\
\textnormal{JSP} & \texttt{la10} & 2 & \textnormal{JSP} & \texttt{la11} & 2 \\
\textnormal{JSP} & \texttt{la12} & 2 & \textnormal{JSP} & \texttt{orb01} & 2 \\
\textnormal{JSP} & \texttt{orb02} & 1 & \textnormal{JSP} & \texttt{orb03} & 2 \\
\textnormal{JSP} & \texttt{orb04} & 2 & \textnormal{JSP} & \texttt{orb05} & 2 \\
\textnormal{MIPLIB} & \texttt{50v-10} & 2 & \textnormal{MIPLIB} & \texttt{CMS750\_4} & 2 \\
\textnormal{MIPLIB} & \texttt{assign1-5-8} & 2 & \textnormal{MIPLIB} & \texttt{atlanta-ip} & 2 \\
\textnormal{MIPLIB} & \texttt{b1c1s1} & 2 & \textnormal{MIPLIB} & \texttt{beasleyC3} & 1 \\
\textnormal{MIPLIB} & \texttt{cbs-cta} & 1 & \textnormal{MIPLIB} & \texttt{chromaticindex512-7} & 2 \\
\textnormal{MIPLIB} & \texttt{cmflsp50-24-8-8} & 3 & \textnormal{MIPLIB} & \texttt{cod105} & 1 \\
\textnormal{MIPLIB} & \texttt{comp21-2idx} & 2 & \textnormal{MIPLIB} & \shortstack[l]{\texttt{cryptanalysis}\\\texttt{kb128n5obj16}} & 3 \\
\textnormal{MIPLIB} & \texttt{dano3\_3} & 1 & \textnormal{MIPLIB} & \texttt{drayage-25-23} & 1 \\
\textnormal{MIPLIB} & \texttt{dws008-01} & 2 & \textnormal{MIPLIB} & \texttt{fhnw-binpack4-48} & 3 \\
\textnormal{MIPLIB} & \shortstack[l]{\texttt{gfd-schedule}\\\texttt{n180f7d50m30k18}} & 3 & \textnormal{MIPLIB} & \texttt{glass-sc} & 2 \\
\textnormal{MIPLIB} & \texttt{glass4} & 2 & \textnormal{MIPLIB} & \texttt{gmu-35-50} & 2 \\
\textnormal{MIPLIB} & \texttt{h80x6320d} & 1 & \textnormal{MIPLIB} & \texttt{irish-electricity} & 3 \\
\textnormal{MIPLIB} & \texttt{lectsched-5-obj} & 2 & \textnormal{MIPLIB} & \texttt{momentum1} & 2 \\
\textnormal{MIPLIB} & \texttt{mushroom-best} & 2 & \textnormal{MIPLIB} & \texttt{mzzv11} & 1 \\
\textnormal{MIPLIB} & \texttt{n3div36} & 2 & \textnormal{MIPLIB} & \texttt{neos-1354092} & 3 \\
\textnormal{MIPLIB} & \texttt{neos-1445765} & 1 & \textnormal{MIPLIB} & \texttt{neos-3046615-murg} & 2 \\
\textnormal{MIPLIB} & \texttt{neos-3083819-nubu} & 1 & \textnormal{MIPLIB} & \texttt{neos-3627168-kasai} & 1 \\
\textnormal{MIPLIB} & \texttt{neos-3988577-wolgan} & 3 & \textnormal{MIPLIB} & \texttt{neos-5052403-cygnet} & 2 \\
\textnormal{MIPLIB} & \texttt{neos-662469} & 2 & \textnormal{MIPLIB} & \texttt{neos-827175} & 1 \\
\textnormal{MIPLIB} & \texttt{neos-911970} & 1 & \textnormal{MIPLIB} & \texttt{neos8} & 1 \\
\textnormal{MIPLIB} & \texttt{neos859080} & 1 & \textnormal{MIPLIB} & \texttt{nexp-150-20-8-5} & 2 \\
\textnormal{MIPLIB} & \texttt{nursesched-sprint02} & 1 & \textnormal{MIPLIB} & \texttt{p200x1188c} & 1 \\
\textnormal{MIPLIB} & \texttt{piperout-27} & 1 & \textnormal{MIPLIB} & \texttt{radiationm18-12-05} & 2 \\
\textnormal{MIPLIB} & \texttt{radiationm40-10-02} & 2 & \textnormal{MIPLIB} & \texttt{ran14x18-disj-8} & 2 \\
\textnormal{MIPLIB} & \texttt{rocI-4-11} & 1 & \textnormal{MIPLIB} & \texttt{rocII-5-11} & 2 \\
\textnormal{MIPLIB} & \texttt{roll3000} & 1 & \textnormal{MIPLIB} & \texttt{s100} & 2 \\
\textnormal{MIPLIB} & \texttt{snp-02-004-104} & 2 & \textnormal{MIPLIB} & \texttt{supportcase22} & 3 \\
\textnormal{MIPLIB} & \texttt{supportcase40} & 2 & \textnormal{MIPLIB} & \texttt{supportcase42} & 2 \\
\textnormal{MIPLIB} & \texttt{swath1} & 1 & \textnormal{MIPLIB} & \texttt{thor50dday} & 2 \\
\textnormal{MIPLIB} & \texttt{traininstance2} & 2 & \textnormal{MIPLIB} & \texttt{uct-subprob} & 2 \\
\textnormal{MIPLIB} & \texttt{unitcal\_7} & 1 & \textnormal{UFL} & \texttt{1023UnifS} & 1 \\
\textnormal{UFL} & \texttt{1123UnifS} & 1 & \textnormal{UFL} & \texttt{1223UnifS} & 1 \\
\textnormal{UFL} & \texttt{123UnifS} & 1 & \textnormal{UFL} & \texttt{1323UnifS} & 1 \\
\textnormal{UFL} & \texttt{1423UnifS} & 1 & \textnormal{UFL} & \texttt{1523UnifS} & 1 \\
\textnormal{UFL} & \texttt{1623UnifS} & 1 & \textnormal{UFL} & \texttt{1723UnifS} & 1 \\
\textnormal{UFL} & \texttt{1823UnifS} & 1 & \textnormal{UFL} & \texttt{1923UnifS} & 1 \\
\textnormal{UFL} & \texttt{2023UnifS} & 1 & \textnormal{UFL} & \texttt{2123UnifS} & 1 \\
\textnormal{UFL} & \texttt{2223UnifS} & 1 & \textnormal{UFL} & \texttt{223UnifS} & 1 \\
\textnormal{UFL} & \texttt{2323UnifS} & 1 & \textnormal{UFL} & \texttt{2423UnifS} & 1 \\
\textnormal{UFL} & \texttt{2523UnifS} & 1 & \textnormal{UFL} & \texttt{2623UnifS} & 1 \\
\textnormal{UFL} & \texttt{2723UnifS} & 1 & \textnormal{UFL} & \texttt{2823UnifS} & 1 \\
\textnormal{UFL} & \texttt{2923UnifS} & 1 & \textnormal{UFL} & \texttt{3023UnifS} & 1 \\
\textnormal{UFL} & \texttt{323UnifS} & 1 & \textnormal{UFL} & \texttt{423UnifS} & 1 \\
\textnormal{UFL} & \texttt{523UnifS} & 1 & \textnormal{UFL} & \texttt{623UnifS} & 1 \\
\textnormal{UFL} & \texttt{723UnifS} & 1 & \textnormal{UFL} & \texttt{823UnifS} & 1 \\
\textnormal{UFL} & \texttt{923UnifS} & 1 &  &  &  \\
\bottomrule
\end{tabular}
}
\floatnotes{Each instance is assigned to a case once under the default configuration $p^{\mathrm{def}}$ with time limit $\tau=600\,\mathrm{s}$. These assignments remain fixed across all experiments.}
\end{table}

\subsection{Within-Family Split Construction}
Exp.~C uses deterministic within-family train/test splits constructed from the Case~1 instances of CVRP, GAP, JSP, and UFL. MIPLIB is excluded because it contains heterogeneous general-purpose instances rather than a single problem family. For each family $s$ and training size $m$, a split is eligible only if at least four Case~1 instances remain for testing. We generate five distinct training subsets for each eligible $(s,m)$ pair using the manifest order and a fixed seed key \texttt{family|m|42}. The corresponding test set consists of all remaining Case~1 instances in the same family. When available, Case~2 and Case~3 instances from that family are also used to evaluate transfer performance, but their logs are never exposed during tuning. Table~\ref{tab:split_design} summarizes the eligible splits and test sets. This design enables within-family transfer comparisons under matched solver budgets while preventing information leakage from test-instance logs into the tuning process.

\begin{table*}[htbp]
\centering
\caption{Deterministic within-family split design for Exp.~C.}
\label{tab:split_design}
\footnotesize
\setlength{\tabcolsep}{4pt}
\renewcommand{\arraystretch}{1.06}
\begin{tabular*}{\textwidth}{@{\extracolsep{\fill}}l r l r l l@{}}
\toprule
\textnormal{Family} & \textnormal{Case~1 pool} & \textnormal{C2 train sizes $m$} & \textnormal{Splits per $m$} & \textnormal{Case~1 test count} & \textnormal{Other transfer targets} \\
\midrule
CVRP & 11 & $2,3,4,5,6$ & 5 & $11-m$ & Case~2 and Case~3 \\
GAP & 7 & $2,3$ & 5 & $7-m$ & Case~2 only \\
JSP & 7 & $2,3$ & 5 & $7-m$ & Case~2 only \\
UFL & 30 & $2,3,4,5,6$ & 5 & $30-m$ & none beyond Case~1 \\
\bottomrule
\end{tabular*}
\floatnotes{All split manifests are generated from the fixed case-assignment manifest using base seed 42. Each eligible split leaves at least four Case~1 instances for testing. C1 uses every Case~1 instance from the listed families as a single anchor.}
\end{table*}

\section{Experimental Protocol, Baselines, and Shared Execution Contract}\label{app_extra_baselines}
This appendix documents the experimental protocol for method comparisons and the controls used in Exp.~B. Table~\ref{tab:execution_contract} specifies the execution conditions shared across the tuning experiments. The following subsections describe the method-specific feedback channels, search-space representations, and prompt modifications. Appendix~\ref{app_prompt_details} provides representative prompts and outputs. Figure~\ref{fig:prompt_example} illustrates the prompt components used by \name{}, Figure~\ref{fig:llm_output_example} presents a representative JSON response, and Table~\ref{tab:prompt_ablation_blocks} identifies the components retained or removed in each Exp.~B control.

\subsection{Shared Execution Contract}
Table~\ref{tab:execution_contract} defines the common execution contract for all reported tuning comparisons, including solver calls, objective accounting, the editable parameter space, advice validation, and randomness controls. These conditions are fixed independently of the method-specific tuning policy.

\begin{table*}[htbp]
\centering
\caption{Shared execution contract for the reported tuning comparisons.}
\label{tab:execution_contract}
\footnotesize
\setlength{\tabcolsep}{4pt}
\renewcommand{\arraystretch}{1.06}
\begin{tabular*}{\textwidth}{@{\extracolsep{\fill}}L{0.26\textwidth} L{0.67\textwidth}@{}}
\toprule
\textnormal{Component} & \textnormal{Contract} \\
\midrule
Solver call & SCIP~9.2.4, one thread, one blocking solver run per evaluated configuration and instance, with per-run time limit $\tau=600\,\mathrm{s}$. \\
Objective and normalization & The regime-adaptive loss in Section~\ref{sec_problem_formulation} is fixed before comparison and used identically by \name{}, LLM controls, and scalar-feedback baselines; each method's own default run (trial~0) defines its normalization reference. \\
Reported selection & For every method, the reported configuration after $t$ trials is the best-so-far configuration over trial~0 (the default) and proposals $1,\dots,t$. \\
Evaluation budget & Each tuning trial consumes one sequential solver evaluation per active instance. Experiment-specific values of $K$ are stated in the corresponding experiment descriptions and split manifests. \\
Editable space & All methods search the same curated whitelist $\mathcal W$ of 45 SCIP parameters in Appendix~\ref{app_whitelist}; numeric entries use the native linear scale. \\
LLM proposal contract & LLM-based methods use edit cap $M=20$, edits around the current best configuration, the same JSON schema, the same re-query limit $V=3$, and the same archive novelty check. \\
Validation & Candidate edits must pass schema validation, whitelist/domain checks, edit-cap truncation, archive-novelty filtering, and successful loading as a SCIP settings file before a solver call. \\
Randomness control & Solver calls fix \texttt{randomization/randomseedshift}=42; stochastic scalar-feedback optimizers use optimizer seed 42. \\
\bottomrule
\end{tabular*}
\end{table*}

\subsection{Black-Box Scalar-Feedback Baselines}
The scalar-feedback baselines receive only the loss returned by the regime-adaptive objective. They have no access to raw solver logs, parsed diagnostics, formulation text, prompt memory, or LLM-side validation feedback. Table~\ref{tab:blackbox_repro} summarizes each method's search-space encoding, feedback channel, and handling of invalid configurations under the shared execution contract in Table~\ref{tab:execution_contract}.

\begin{table*}[htbp]
\centering
\caption{Implementation details of the scalar-feedback baselines.}
\label{tab:blackbox_repro}
\footnotesize
\setlength{\tabcolsep}{3pt}
\renewcommand{\arraystretch}{1.08}
\begin{tabular*}{\textwidth}{@{\extracolsep{\fill}}l L{0.30\textwidth} c L{0.35\textwidth}@{}}
\toprule
\textnormal{Method} & \textnormal{Search-space encoding} & \textnormal{Feedback} & \textnormal{Settings and invalid handling} \\
\midrule
\textsc{SMAC3} & Typed \textsc{ConfigSpace} over the whitelist; native linear numeric domains. & scalar loss & Deterministic, seed 42, $K$ sequential evaluations; non-finite values noncompetitive. \\
\textsc{RBFOpt} & Native float/int bounds; Boolean and categorical dimensions integer-encoded. & scalar loss & Seed 42, no local search, genetic global search, $K$ sequential evaluations; non-finite values noncompetitive. \\
\textsc{Random} & Independent typed sampling; positive real intervals log-uniform, all other numeric intervals uniform. & scalar loss & Seed 42, $K$ sequential samples, no surrogate; the default enters only through trial~0. \\
\bottomrule
\end{tabular*}
\floatnotes{All three methods optimize only the scalar loss returned by the current case objective.}
\end{table*}

\subsection{Prompt-Based LLM Controls}
\begin{figure*}[htbp]
\centering
\caption{Archived prompt excerpt for CVRP instance A-n37-k6 after trial~2.}
\label{fig:prompt_trial2_screenshot}
\includegraphics[width=0.99\textwidth]{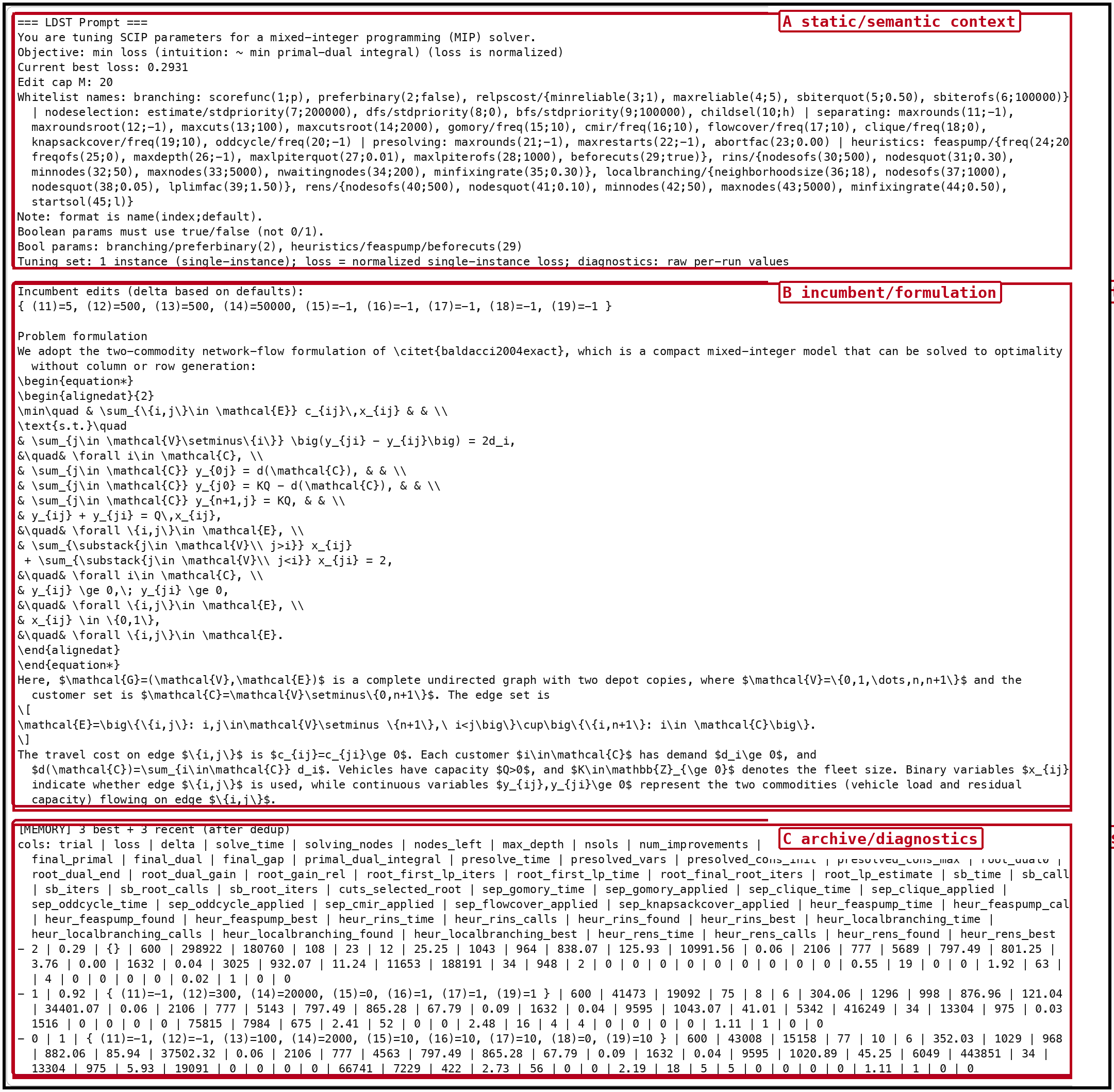}
\floatnotes{The cropped screenshot is taken from the archived prompt and spans the prompt header through the archive-memory block. Red boxes identify the semantic and static context (A), the incumbent (current best configuration) and formulation context (B), and archive memory containing parsed diagnostics (C).}
\end{figure*}

Let $\Pi^{\textsc{LDST}}$ denote the full \name{} prompt. It comprises (i) a static header specifying the objective, parameter whitelist, current best configuration, edit cap, and output schema, (ii) optional formulation text, (iii) validator feedback during re-queries, (iv) archive memory, (v) a deduplication list recording previously evaluated configurations, and (vi) a diagnostic-feature guide linking parsed log fields to potential solver bottlenecks. Figure~\ref{fig:prompt_trial2_screenshot} presents a cropped prompt screenshot from an actual Exp.~A run, highlighting the components relevant to the Exp.~B controls. The three LLM controls selectively remove or restrict information available to the decision policy while preserving the same outer tuning loop.

\paragraph{LLM-StaticContext.}
\textsc{LLM-StaticContext} retains the static components of $\Pi^{\textsc{LDST}}$ and the parsed diagnostics of the current best run. These include whitelist parameter names, defaults, and domains, edits defining the current best configuration, task and edit constraints, optional formulation text, and a single block of current diagnostics. In Figure~\ref{fig:prompt_trial2_screenshot}, this control retains regions A and B but removes the cross-trial archive in C, the diagnostic guide, and the deduplication list. Thus, the LLM receives information about the current best configuration and its solver behavior without access to the history of previous trials.

\paragraph{LLM-OutcomeOnly.}
\textsc{LLM-OutcomeOnly} replaces the full archive in $\Pi^{\textsc{LDST}}$ with a scalar-only memory $\mathcal M^{\textsc{scalar}}$ containing the trial id, loss, solve time, final gap, and finite/feasible status whenever these quantities are defined. In Figure~\ref{fig:prompt_trial2_screenshot}, this control retains region A and the scalar feedback in region C but removes the parsed diagnostic fields and the diagnostic guide. It therefore preserves outcome feedback across trials while removing process-level log structure.

\paragraph{LLM-IndexOnly.}
\textsc{LLM-IndexOnly} retains the archive, log-derived diagnostics, diagnostic guide, and deduplication list but restricts the semantic parameter context. Relative to Figure~\ref{fig:prompt_trial2_screenshot}, it retains region C while replacing the parameter names, groups, and defaults in region A with editable parameter indices, per-index domain hints, and the indices corresponding to Boolean parameters. These hints support valid parameter edits without revealing their meanings. The tuning loop, validator, and JSON schema remain unchanged, except that the optional textual \texttt{whitelist} field is omitted from the returned edits. In Exp.~B, \textsc{LLM-IndexOnly} terminates early on \texttt{gap\_gapc-4}, \texttt{gap\_gapd-1}, \texttt{jsp\_la01}, \texttt{miplib\_neos-827175}, \texttt{ufl\_1023UnifS}, and \texttt{ufl\_2423UnifS} after trials 20, 16, 18, 23, 20, and 21, respectively. \textsc{LLM-StaticContext} also terminates early on \texttt{miplib\_cbs-cta} after trial~26.

\paragraph{Formulation and leakage control.}
All reported LLM-based runs use \texttt{formulation=auto}, which includes the built-in formulation text for CVRP, GAP, JSP, and UFL and none for MIPLIB. The formulation text is a family-level mathematical description retrieved from \texttt{utilities/formulations.py}. It depends only on the benchmark family and is identical for all instances within that family. No formulation text is constructed using test-instance logs or target-instance solve trajectories. In Exp.~C, the LLM receives only logs from the training instances during tuning. Test instances from Cases~1, 2, and~3 are evaluated only after a configuration has been selected, and their logs are never fed back into the tuning prompts. Table~\ref{tab:method_info_matrix} summarizes the information available to the black-box baselines, LLM controls, and full \name{} framework.

\begin{table*}[htbp]
\centering
\caption{Information-access matrix for the main control families.}
\label{tab:method_info_matrix}
\scriptsize
\setlength{\tabcolsep}{2pt}
\renewcommand{\arraystretch}{1.05}
\begin{tabular*}{\textwidth}{@{\extracolsep{\fill}}lcccccc@{}}
\toprule
\textnormal{Method} & \textnormal{Scalar hist.} & \textnormal{Parsed diag.} & \textnormal{Param. sem.} & \textnormal{Form. text} & \textnormal{Rationale} & \textnormal{Validator} \\
\midrule
\textnormal{BBO baselines} & \checkmark & \xmark & \xmark & \xmark & \xmark & \xmark \\
\textsc{LLM-StaticContext} & \xmark & \checkmark & \checkmark & optional & \checkmark & \checkmark \\
\textsc{LLM-OutcomeOnly} & \checkmark & \xmark & \checkmark & optional & \checkmark & \checkmark \\
\textsc{LLM-IndexOnly} & \checkmark & \checkmark & \xmark & optional & \checkmark & \checkmark \\
\textsc{\name{}} & \checkmark & \checkmark & \checkmark & optional & \checkmark & \checkmark \\
\bottomrule
\end{tabular*}
\floatnotes{The BBO row groups \textsc{SMAC3}, \textsc{RBFOpt}, and \textsc{Random}. Column headers abbreviate scalar history, parsed diagnostics, parameter semantics, optional formulation text, visible rationale field, and the prompt-side validator, re-query, and deduplication contract. Black-box baselines sample directly from the typed whitelist and do not use this prompt-side contract. \textsc{SMAC3} and \textsc{RBFOpt} fit surrogate models to scalar history, whereas \textsc{Random} ignores that history. For \textsc{LLM-IndexOnly}, index-level domain hints are not counted as parameter semantics because parameter names and defaults are hidden.}
\end{table*}

\section{Configuration Space: SCIP Parameter Whitelist}\label{app_whitelist}
We use the same curated whitelist $\mathcal{W}$ across all experiments and controls. The whitelist is small enough that edits remain interpretable and semantically checkable, while covering branching, node selection, cut separation, presolving, and primal heuristics, the main solver components reported in SCIP logs. Table~\ref{tab:whitelist_full} provides the complete specification of all 45 parameters used in the implementation, including their SCIP parameter paths, types, default values, allowed domains, and scales. All numeric parameters are edited on their native linear scales, without parameter-specific transformations. In addition to the typed whitelist, the LLM-side validator checks JSON-schema validity, truncates proposals to the edit cap, checks archive novelty, and verifies that the proposed configuration loads as a SCIP settings file.

\begin{table}[htbp]
\centering
\caption{Full typed SCIP parameter whitelist used in the experiments.}
\label{tab:whitelist_full}
{\scriptsize
\setlength{\tabcolsep}{2pt}
\renewcommand{\arraystretch}{1.0}
\begin{tabular*}{\textwidth}{@{\extracolsep{\fill}}L{0.40\textwidth} c c C{0.30\textwidth} c@{}}
\toprule
\textnormal{SCIP parameter path} & \textnormal{Type} & \textnormal{Default} & \textnormal{Allowed domain} & \textnormal{Scale} \\
\midrule
\path|branching/scorefunc| & cat & \texttt{p} & \texttt{\{s, p, q\}} & lin. \\
\path|branching/preferbinary| & bool & \texttt{false} & \texttt{\{false, true\}} & lin. \\
\path|branching/relpscost/minreliable| & float & \texttt{1} & \texttt{[0, 1000000000]} & lin. \\
\path|branching/relpscost/maxreliable| & float & \texttt{5} & \texttt{[0, 1000000000]} & lin. \\
\path|branching/relpscost/sbiterquot| & float & \texttt{0.5} & \texttt{[0, 1000000000]} & lin. \\
\path|branching/relpscost/sbiterofs| & int & \texttt{100000} & \texttt{[0, 1000000000]} & lin. \\
\path|nodeselection/estimate/stdpriority| & int & \texttt{200000} & \texttt{[-536870912, 1000000000]} & lin. \\
\path|nodeselection/dfs/stdpriority| & int & \texttt{0} & \texttt{[-536870912, 1000000000]} & lin. \\
\path|nodeselection/bfs/stdpriority| & int & \texttt{100000} & \texttt{[-536870912, 1000000000]} & lin. \\
\path|nodeselection/childsel| & cat & \texttt{h} & \texttt{\{d, u, p, i, l, r, h\}} & lin. \\
\path|separating/maxrounds| & int & \texttt{-1} & \texttt{[-1, 1000000000]} & lin. \\
\path|separating/maxroundsroot| & int & \texttt{-1} & \texttt{[-1, 1000000000]} & lin. \\
\path|separating/maxcuts| & int & \texttt{100} & \texttt{[0, 1000000000]} & lin. \\
\path|separating/maxcutsroot| & int & \texttt{2000} & \texttt{[0, 1000000000]} & lin. \\
\path|separating/gomory/freq| & int & \texttt{10} & \texttt{[-1, 1000000000]} & lin. \\
\path|separating/cmir/freq| & int & \texttt{10} & \texttt{[-1, 1000000000]} & lin. \\
\path|separating/flowcover/freq| & int & \texttt{10} & \texttt{[-1, 1000000000]} & lin. \\
\path|separating/clique/freq| & int & \texttt{0} & \texttt{[-1, 1000000000]} & lin. \\
\path|separating/knapsackcover/freq| & int & \texttt{10} & \texttt{[-1, 1000000000]} & lin. \\
\path|separating/oddcycle/freq| & int & \texttt{-1} & \texttt{[-1, 1000000000]} & lin. \\
\path|presolving/maxrounds| & int & \texttt{-1} & \texttt{[-1, 1000000000]} & lin. \\
\path|presolving/maxrestarts| & int & \texttt{-1} & \texttt{[-1, 1000000000]} & lin. \\
\path|presolving/abortfac| & float & \texttt{0.0008} & \texttt{[0, 1]} & lin. \\
\path|heuristics/feaspump/freq| & int & \texttt{20} & \texttt{[-1, 1000000000]} & lin. \\
\path|heuristics/feaspump/freqofs| & int & \texttt{0} & \texttt{[0, 1000000000]} & lin. \\
\path|heuristics/feaspump/maxdepth| & int & \texttt{-1} & \texttt{[-1, 1000000000]} & lin. \\
\path|heuristics/feaspump/maxlpiterquot| & float & \texttt{0.01} & \texttt{[0, 1000000000]} & lin. \\
\path|heuristics/feaspump/maxlpiterofs| & int & \texttt{1000} & \texttt{[0, 1000000000]} & lin. \\
\path|heuristics/feaspump/beforecuts| & bool & \texttt{true} & \texttt{\{false, true\}} & lin. \\
\path|heuristics/rins/nodesofs| & int & \texttt{500} & \texttt{[0, 1000000000]} & lin. \\
\path|heuristics/rins/nodesquot| & float & \texttt{0.3} & \texttt{[0, 1]} & lin. \\
\path|heuristics/rins/minnodes| & int & \texttt{50} & \texttt{[0, 1000000000]} & lin. \\
\path|heuristics/rins/maxnodes| & int & \texttt{5000} & \texttt{[0, 1000000000]} & lin. \\
\path|heuristics/rins/nwaitingnodes| & int & \texttt{200} & \texttt{[0, 1000000000]} & lin. \\
\path|heuristics/rins/minfixingrate| & float & \texttt{0.3} & \texttt{[0, 1]} & lin. \\
\path|heuristics/localbranching/neighborhoodsize| & int & \texttt{18} & \texttt{[1, 1000000000]} & lin. \\
\path|heuristics/localbranching/nodesofs| & int & \texttt{1000} & \texttt{[0, 1000000000]} & lin. \\
\path|heuristics/localbranching/nodesquot| & float & \texttt{0.05} & \texttt{[0, 1]} & lin. \\
\path|heuristics/localbranching/lplimfac| & float & \texttt{1.5} & \texttt{[1, 1000000000]} & lin. \\
\path|heuristics/rens/nodesofs| & int & \texttt{500} & \texttt{[0, 1000000000]} & lin. \\
\path|heuristics/rens/nodesquot| & float & \texttt{0.1} & \texttt{[0, 1]} & lin. \\
\path|heuristics/rens/minnodes| & int & \texttt{50} & \texttt{[0, 1000000000]} & lin. \\
\path|heuristics/rens/maxnodes| & int & \texttt{5000} & \texttt{[0, 1000000000]} & lin. \\
\path|heuristics/rens/minfixingrate| & float & \texttt{0.5} & \texttt{[0, 1]} & lin. \\
\path|heuristics/rens/startsol| & cat & \texttt{l} & \texttt{\{n, l\}} & lin. \\
\bottomrule
\end{tabular*}
}
\floatnotes{The whitelist is generated directly from the code under SCIP~9.2.4.}
\end{table}

\section{LLM Hyperparameters and Shared Settings}\label{app_hparams}
The LLM serves only as a constrained proposal mechanism. Table~\ref{tab:exp_hparams} summarizes the model-side settings not specified in the shared execution contract in Appendix~\ref{app_extra_baselines}. Experiment-specific trial counts are provided in the corresponding experiment descriptions and Appendix~\ref{app_artifact}.
\begin{table}[htbp]
\centering
\caption{Model-side hyperparameters for the LLM-based tuning loops.}
\label{tab:exp_hparams}
\footnotesize
\setlength{\tabcolsep}{4pt}
\renewcommand{\arraystretch}{1.05}
\begin{tabular}{@{}p{0.52\linewidth} p{0.40\linewidth}@{}}
\toprule
Item & Value \\
\midrule
Best and recent records in prompt memory $(B,R)$ & 3, 3 \\
LLM model & GPT-5.2 \\
Reasoning effort & medium \\
Max output tokens & 20{,}000 \\
\bottomrule
\end{tabular}
\end{table}

\section{SCIP Log Parser and Diagnostics}\label{app_scip_diagnostics}
For each run $k$, the parser converts the raw SCIP log $\Lambda^{(k)}$ into a diagnostic vector with a fixed schema,
\[
\mathbf d^{(k)}=\operatorname{Parse}(\Lambda^{(k)}),
\]
which is passed to the prompt builder described in Section~\ref{sec_methodology}. The parser requires no access to internal solver states and uses only standard SCIP log sections, including the branch-and-bound progress table, root-node report, presolve summary, heuristic table, separator table, and final solving summary. Figure~\ref{fig:log_parse_bridge} illustrates how the raw SCIP log is converted into the compact diagnostic state used by \name{}.

\begin{figure*}[htbp]
  \centering
  \caption{Representative raw SCIP log excerpt and corresponding parsed diagnostic summary.}
  \label{fig:log_parse_bridge}
  \begin{minipage}{0.485\textwidth}
  \footnotesize
  \begin{tcolorbox}[
    colback=white,
    colframe=black!70,
    boxrule=0.6pt,
    arc=1.2mm,
    left=1.5mm,right=1.5mm,top=1.0mm,bottom=1.0mm
  ]
\begin{verbatim}
=== RAW SCIP LOG (excerpt) ===
presolved problem has 953 vars, 358 cons
presolving time: 0.12

 time | node  | left  | gap
 0.2s |     1 |     0 | 12.44%
 0.3s |     1 |     0 |  4.50%
 7.5s |  4200 |  3971 |  0.23%
10.1s |  6691 |  6318 |  0.17%
23.7s | 18580 | 17424 |  0.11%
\end{verbatim}
  \end{tcolorbox}
  \end{minipage}\hfill
  \begin{minipage}{0.485\textwidth}
  \footnotesize
  \begin{tcolorbox}[
    colback=white,
    colframe=black!70,
    boxrule=0.6pt,
    arc=1.2mm,
    left=1.5mm,right=1.5mm,top=1.0mm,bottom=1.0mm
  ]
\begin{verbatim}
=== PARSED DIAGNOSTICS (excerpt) ===
t_total            = 600.0
nodes_last         = 432241
nodes_left         = 319848
max_depth          = 108
num_improvements   = 12
first_inc_time     = 0.0
gap_final          = 1500.39
pdi                = 76.83
presolve_time      = 0.12
root_dual_gain     = 9.06e-03
sb_time            = 6.05
sb_iters           = 116575
cuts_selected_total= 95704
sep_cmir_applied   = 29937
heur_rins_time     = 0.41
\end{verbatim}
  \end{tcolorbox}
  \end{minipage}
\floatnotes{The left box shows an excerpt from a real Case~2 run. The parser also reads the corresponding bound, LP, cut, separator, heuristic, and summary sections from the same log. The right box shows the compact, fixed-schema diagnostic state returned by the parser and consumed by \name{}.}
\end{figure*}

Table~\ref{tab:diag_schema} defines the diagnostic quantities used by the prompt builder. Let $\mathcal R^{(k)}$ denote the ordered progress-table rows parsed from run $\Lambda^{(k)}$. Each row $r\in\mathcal R^{(k)}$ contains elapsed time $t_r$, processed-node count $n_r$, primal bound $z^{\mathrm{pr}}_r$, and dual bound $z^{\mathrm{du}}_r$ when these fields are available. For any row field $x$, let $\operatorname{first}_{f}(x)$ and $\operatorname{last}_{f}(x)$ denote its first and last finite parsed values, respectively, and set $\eta=10^{-12}$. Let $\mathcal S=\{\mathrm{gomory},\mathrm{clique},\mathrm{oddcycle},\mathrm{cmir},\mathrm{flowcover},\mathrm{knapsackcover}\}$ and $\mathcal H=\{\mathrm{feaspump},\mathrm{rins},\mathrm{localbranching},\mathrm{rens}\}$ denote the tracked separator and heuristic families, respectively. Summary quantities are extracted from the final SCIP summary when available, with the row-based fallback rules specified in the table.
\begin{table}[htbp]
\centering
\caption{Diagnostic quantities used by the prompt builder.}
\label{tab:diag_schema}
\scriptsize
\setlength{\tabcolsep}{3pt}
\renewcommand{\arraystretch}{1.08}
\begin{tabular*}{\textwidth}{@{\extracolsep{\fill}}L{0.15\textwidth} L{0.28\textwidth} L{0.39\textwidth} L{0.13\textwidth}@{}}
\toprule
\textnormal{Block} & \textnormal{Parser fields exposed to the policy} & \textnormal{Implementation definition} & \textnormal{SCIP source} \\
\midrule
Run size and status &
\texttt{status}, \texttt{t\_total} (shown as \texttt{solve\_time}), \texttt{nodes\_last} (shown as \texttt{solving\_nodes}), \texttt{nodes\_left}, \texttt{max\_depth} &
\texttt{t\_total} is the reported solving time, or $\operatorname{last}_{f}(t)$ if missing; \texttt{nodes\_last} is the reported solving-node count, or $\operatorname{last}_{f}(n)$ if missing. \texttt{nodes\_left} and \texttt{max\_depth} are the B\&B-tree summary values. &
final summary, B\&B tree, progress rows \\
\midrule
Bounds and integral &
\texttt{primal\_final}, \texttt{dual\_final}, \texttt{gap\_final}, \texttt{pdi} (shown as \texttt{primal\_dual\_integral}) &
\texttt{primal\_final} and \texttt{dual\_final} are final-summary bounds, falling back to $\operatorname{last}_{f}(z^{\mathrm{pr}})$ and $\operatorname{last}_{f}(z^{\mathrm{du}})$. \texttt{gap\_final}$=|\texttt{primal\_final}-\texttt{dual\_final}|$. \texttt{pdi}$=\mathrm{PDI}^{\mathrm{SCIP}}(\Lambda^{(k)})$ is the SCIP-reported primal--dual integral parsed from the Integrals block. &
final summary, Integrals block \\
\midrule
Incumbents &
\texttt{nsols}, \texttt{num\_improvements}, \texttt{first\_inc\_time} (shown as \texttt{first\_incumbent\_time}), \texttt{first\_primal} &
\texttt{nsols} and \texttt{num\_improvements} are solution-summary counts. \texttt{first\_inc\_time} and \texttt{first\_primal} use the first-solution summary when present; otherwise they use the time and value at $\operatorname{first}_{f}(z^{\mathrm{pr}})$. &
solution summary, progress rows \\
\midrule
Presolve &
\texttt{presolve\_time}, \texttt{presolved\_vars}, \texttt{presolved\_cons\_init}, \texttt{presolved\_cons\_max} &
Direct parser entries for presolving time, presolved-variable count, initial presolved-constraint count, and maximum presolved-constraint count. &
presolve summary \\
\midrule
Root progress &
\texttt{root\_dual0}, \texttt{root\_dual\_end}, \texttt{root\_dual\_gain}, \texttt{root\_gain\_rel}, \texttt{root\_first\_lp\_iters}, \texttt{root\_first\_lp\_time}, \texttt{root\_final\_root\_iters}, \texttt{root\_lp\_estimate} &
\texttt{root\_dual0} is the first root-LP value and \texttt{root\_dual\_end} is the final root dual bound. The derived fields are \texttt{root\_dual\_gain}$=\texttt{root\_dual\_end}-\texttt{root\_dual0}$ and \texttt{root\_gain\_rel}$=\texttt{root\_dual\_gain}/\max\{|\texttt{root\_dual0}|,\eta\}$. &
root-node report \\
\midrule
LP and branching effort &
\texttt{sb\_time}, \texttt{sb\_calls}, \texttt{sb\_iters}, \texttt{sb\_root\_calls}, \texttt{sb\_root\_iters} &
\texttt{sb\_time}, \texttt{sb\_calls}, and \texttt{sb\_iters} are parsed from the LP-statistics row for strong branching; \texttt{sb\_root\_calls} and \texttt{sb\_root\_iters} are parsed from the root-node LP row. &
LP statistics \\
\midrule
Cut separation &
\texttt{cuts\_selected\_total}, \texttt{cuts\_selected\_root}, and applied-cut or time fields for $c\in\mathcal S$ &
\texttt{cuts\_selected\_total}$=\sum_q \texttt{selected}_q$ and \texttt{cuts\_selected\_root}$=\sum_q \texttt{rootselec}_q$ over cut selectors $q$. For each separator family $c$, the parser records the applied-cut count and records execution time when the SCIP table exposes that entry. &
cut-selector and separator tables \\
\midrule
Heuristic activity &
time, call, found-solution, and best-solution counts for $h\in\mathcal H$ &
For each heuristic family $h$, the parser exposes SCIP's \texttt{calls}, \texttt{exectime}, \texttt{found}, and \texttt{best} columns as prompt diagnostics. &
primal-heuristics table \\
\bottomrule
\end{tabular*}
\floatnotes{The quantities are parsed from standard SCIP logs. Aliases in parentheses give the field names rendered in the archived prompts.}
\end{table}

\noindent The parser converts detailed solver logs into a compact diagnostic state for the decision policy. In single-instance mode, the resulting summary is included directly in the prompt archive. In multi-instance mode, the same schema is applied to each instance, and the resulting diagnostics are presented as aligned rows. This preserves comparability across instances without requiring the LLM to process raw solver logs.
\section{Prompt Construction and Representative Outputs}\label{app_prompt_details}
This section presents a representative prompt and LLM response and explains their relationship to the prompt ablations in Exp.~B. It specifies the information available to the decision policy, the components removed in each ablation, and how the proposed parameter edits are supported by log-derived evidence.

\begin{figure}[htbp]
  \centering
\caption{Representative single-instance prompt excerpt issued after trial~2.}
  \label{fig:prompt_example}
  \begin{minipage}{0.95\textwidth}
  \footnotesize
  \begin{tcolorbox}[
    colback=white,
    colframe=black!70,
    boxrule=0.6pt,
    arc=1.2mm,
    left=1.5mm,right=1.5mm,top=1.0mm,bottom=1.0mm
  ]
\begin{verbatim}
=== LDST Prompt ===
You are tuning SCIP parameters for a mixed-integer programming solver.
Objective: min loss (intuition: ~ min solution time for opt) (loss is normalized)
Current best loss: 0.4436
Edit cap M: 20
Whitelist names: branching: scorefunc(1;p), preferbinary(2;false), ...
| nodeselection: ... | separating: ... | presolving: ... | heuristics: ...
Tuning set: 1 instance (single-instance);
loss = normalized single-instance loss; diagnostics: raw per-run values

Incumbent edits (delta based on defaults):
{ (2)=true, (4)=1, (5)=0.10, (6)=5000, (11)=0, (12)=50,
  (13)=200, (14)=5000, ... }

Problem formulation
We adopt the two-commodity network-flow formulation of \citet{baldacci2004exact} ...

[MEMORY] 3 best + 3 recent (after dedup)
cols: trial | loss | delta | solve_time | solving_nodes | nodes_left | ...
- 2 | 0.44 | {} | 1.18 | 1074 | 0 | 26 | ...
- 1 | 0.57 | { (4)=2, (5)=0.20, (6)=20000, (11)=-1, (12)=-1 }
  | 1.52 | ...
- 0 | 1 | { (2)=false, (4)=5, (5)=0.50, (6)=100000, ... }
  | 2.66 | ...

[TIPS]
Use the diagnostics and memory to identify ONE dominant bottleneck before ...

[DEDUP]
Avoid these (loss, delta based on incumbent)s other than already listed above
(most recent first):
(empty)

[TASK]
Propose the next configuration as a sparse edit set Delta p relative to ...
Use <= 20 edits and only whitelist indices.
Return ONLY valid JSON with this schema:
{"edits": [{"index": "<int>", "whitelist": "<str>", "value": "<val>"}],
 "rationale": "<brief reasoning tied to the observed diagnostics>"}
\end{verbatim}
  \end{tcolorbox}
  \end{minipage}
\floatnotes{The excerpt is abridged at ``...''. The prompt includes static context, edits defining the current best configuration (labeled ``incumbent'' in the prompt), formulation text when available, archive memory, the diagnostic-feature guide, the deduplication list, and the output schema.}
\end{figure}

\noindent Figure~\ref{fig:prompt_example} illustrates how the prompt is constructed around the current best configuration rather than the most recent trial. The archive summarizes both the best and the most recent runs using a fixed schema, providing the LLM with compact, comparable information across trials instead of raw solver logs.

\begin{table}[htbp]
\centering
\caption{Prompt-content ablations used in Exp.~B.}
\label{tab:prompt_ablation_blocks}
\scriptsize
\setlength{\tabcolsep}{3pt}
\renewcommand{\arraystretch}{1.08}
\begin{tabular*}{\textwidth}{@{\extracolsep{\fill}}L{0.31\textwidth} L{0.15\textwidth} L{0.17\textwidth} L{0.15\textwidth} L{0.15\textwidth}@{}}
\toprule
\textnormal{Prompt block in Figure~\ref{fig:prompt_example}} & \textnormal{\name{}} & \textnormal{\textsc{StaticContext}} & \textnormal{\textsc{OutcomeOnly}} & \textnormal{\textsc{IndexOnly}} \\
\midrule
\textnormal{Header, whitelist, current best, edit cap, schema} & semantic names, defaults, domains & same & same & indices and domains \\
\textnormal{Optional formulation text} & optional & optional & optional & optional \\
\textnormal{Validator feedback during re-query} & kept & re-query only & kept & kept \\
\texttt{[MEMORY]} archive summary & scalar + diagnostics & current diagnostics only & scalar only & full, indexed \\
\textnormal{Multi-instance diagnostic rows} & kept & removed & removed & kept \\
\texttt{[TIPS]} diagnostic-feature guide & kept & removed & removed & kept \\
\texttt{[DEDUP]} prior configurations & kept & removed & kept & kept \\
\bottomrule
\end{tabular*}
\floatnotes{The controls are illustrated using the same prompt excerpt as Figure~\ref{fig:prompt_example}. Each control retains the same outer tuning loop and modifies only the specified prompt blocks. The ``current diagnostics only'' control removes archive memory while retaining the diagnostic block of the current best run. The ``scalar only'' control retains trial id, loss, solve time, final gap, and finite/feasible status. The ``full, indexed'' control retains the diagnostic archive but replaces semantic parameter names with indices.}
\end{table}

\noindent Table~\ref{tab:prompt_ablation_blocks} summarizes the information channels examined in Exp.~B. All variants share the same tuning loop and differ only in their access to process diagnostics, scalar history, and parameter semantics.

\begin{figure}[htbp]
  \centering
\caption{Representative LLM response with a sparse parameter edit set.}
  \label{fig:llm_output_example}
  \begin{minipage}{0.95\textwidth}
  \footnotesize
  \begin{tcolorbox}[
    colback=white,
    colframe=black!70,
    boxrule=0.6pt,
    arc=1.2mm,
    left=1.5mm,right=1.5mm,top=1.0mm,bottom=1.0mm
  ]
\begin{verbatim}
=== LLM_RESPONSE_JSON (excerpt) ===
{
  "edits": [
    {"index":"3", "whitelist":"branching: relpscost/minreliable", "value":"2"},
    {"index":"4", "whitelist":"branching: relpscost/maxreliable", "value":"4"},
    {"index":"5", "whitelist":"branching: relpscost/sbiterquot", "value":"0.05"},
    {"index":"6", "whitelist":"branching: relpscost/sbiterofs", "value":"50"}
  ],
  "rationale":
  "1) Incumbent already fixes the dominant issue at  root
   (root_dual_gain 4.84, nodes drop to 126), and separation times are tiny;
   so root cutting is not the bottleneck anymore.
   2) The remaining dominant cost is expensive branching:
   sb_time=16 out of 28 s with very large sb_iters (128118) despite the small tree.
   3) Therefore, reduce strong-branching effort by making pseudocosts
   become 'reliable' sooner (lower minreliable/maxreliable) and by
   tightening the SB LP-iteration limits (lower sbiterquot/sbiterofs).
   This should cut sb_time substantially
   while accepting a modest risk of more nodes."
}
\end{verbatim}
  \end{tcolorbox}
  \end{minipage}
\floatnotes{The response comes from a different single-instance run than Figure~\ref{fig:prompt_example}. The proposed edits are tied to the parsed diagnostics of that run, making the proposal interpretable and allowing it to be validated before the next solver call.}
\end{figure}

\section{Validation of the Case~1 Completion Signal}\label{app_extrapolation}
This appendix describes the Case~1 completion-time estimator $\widehat{T}_{\infty}$ as implemented and evaluates its absolute accuracy and ranking fidelity. All quantities are derived solely from standard SCIP logs, without access to internal solver states or modifications to the solver source code. In both retrospective studies, we simulate truncation of a completed long run at a specified cutoff and provide the estimator with only the progress-table rows and primal and dual bounds available at that time.

Figure~\ref{fig:extrapolation_log_excerpt} presents an XML100 SCIP log excerpt used in the retrospective validation. Although the full run proves optimality, it still has a nonzero gap at the simulated cutoff $\tau=600\,\mathrm{s}$. The estimator therefore uses only the progress-table fields available up to that point, including elapsed time, processed and open node counts, and primal and dual bounds. The final summary of the completed run provides the true completion time against which the estimate is evaluated.

\begin{figure}[htbp]
\centering
\caption{SCIP log excerpt for XML100 instance \texttt{XML100\_3311\_15\_35}.}
\label{fig:extrapolation_log_excerpt}
\scriptsize
\begin{minipage}{0.96\linewidth}
\begin{tcolorbox}[colback=white,colframe=black,boxrule=0.45pt,arc=2pt,left=4pt,right=4pt,top=4pt,bottom=4pt]
\ttfamily
SCIP> read XML100\_3311\_15\_35.mps.gz

\medskip
\noindent
\setlength{\fboxsep}{3pt}
\begin{tabular}{@{}c@{\hspace{8pt}}c@{}}
\fcolorbox{red!75!black}{white}{%
\begin{minipage}{0.31\linewidth}
\ttfamily
{\normalfont\color{red!75!black} snapshot fields $(t_j,n_j,u_j)$\par}
\ttfamily
\begin{tabular}{@{}rrr@{}}
time & node & left\\
 420s & 41100 & 1305\\
 480s & 47300 & 1237\\
 540s & 52900 & 1069\\
 580s & 56400 &  883\\
 600s & 57900 &  725
\end{tabular}
\end{minipage}}&
\fcolorbox{red!75!black}{white}{%
\begin{minipage}{0.52\linewidth}
\ttfamily
{\normalfont\color{red!75!black} bound fields at each snapshot\par}
\ttfamily
\begin{tabular}{@{}rrrr@{}}
dualbound & primalbound & gap & compl.\\
1.629251e+04 & 1.645100e+04 & 0.97\% & 49.99\%\\
1.630569e+04 & 1.645100e+04 & 0.89\% & 56.04\%\\
1.632154e+04 & 1.645100e+04 & 0.79\% & 67.09\%\\
1.633211e+04 & 1.645100e+04 & 0.73\% & 73.33\%\\
1.633899e+04 & 1.645100e+04 & 0.69\% & 77.45\%
\end{tabular}
\end{minipage}}
\end{tabular}

\medskip
\noindent
\fcolorbox{red!75!black}{white}{%
\begin{minipage}{0.92\linewidth}
\ttfamily
{\normalfont\color{red!75!black} final summary used to audit completion and final primal value\par}
\ttfamily
SCIP Status\hspace{2.8em}: problem is solved [optimal solution found]\par
Solving Time (sec) : 636.04\par
Primal Bound\hspace{2.7em}: +1.64510000000000e+04\par
Dual Bound\hspace{3.4em}: +1.64510000000000e+04\par
Gap\hspace{7.1em}: 0.00 \%
\end{minipage}}
\end{tcolorbox}
\end{minipage}
\floatnotes{The completed run proves optimality at $636.04\,\mathrm{s}$, but its simulated truncation at $\tau=600\,\mathrm{s}$ retains a $0.69\%$ gap and $725$ open nodes. Red boxes identify the progress-table fields used to construct snapshots and the final summary used to record the true completion time.}
\end{figure}

\subsection{Completion-Time Estimator}
The estimator is implemented as \texttt{predict\_time()} in \texttt{utilities/est\_time.py}. If a run is solved within the time limit $\tau$, the estimator returns the observed solve time. Otherwise, it truncates the parsed SCIP progress table at $\tau'=\min\{\tau,t_{\max}\}$, where $t_{\max}$ is the last logged time. After discarding the startup prefix with \texttt{node}$\le 1$, it retains only rows with finite values for time, processed-node count, open-node count, primal bound, and dual bound. The resulting snapshots are $\mathcal S=\{(t_j,n_j,u_j,g_j)\}_{j=1}^m$, where $t_j$ is elapsed time, $n_j$ is the processed-node count, $u_j$ is the number of open nodes, and $g_j=\lvert z_{\mathrm{final}}^{\mathrm{pr}}-z^{\mathrm{du}}(t_j)\rvert$ is the absolute difference between the primal bound reported at the end of the observed run and the dual bound at time $t_j$. For tuning runs, the observed run ends at $\tau$, so $z_{\mathrm{final}}^{\mathrm{pr}}$ is the primal bound available at the time limit.

\begingroup
\predisplaypenalty=100
Following the single-variable branching (SVB) model of \citet{le2017abstract}, each snapshot yields a candidate terminal tree size $\widehat N_j(x)=n_j+u_jx^{g_j}$. The implementation draws up to $25$ stratified samples from the valid rows, requiring at least $8$ rows for fitting. For a sampled subset $\mathcal T\subseteq\mathcal S$, it solves the following optimization problem on the log scale:
\[
\begin{aligned}
\min_{x,b,\{z_j\}_{j\in\mathcal T},v}\quad
& \sum_{j\in\mathcal T} z_j + \lambda v \\
\text{s.t.}\quad
& z_j \ge \log b - \log \widehat N_j(x), && \forall j\in\mathcal T,\\
& z_j \ge \log \widehat N_j(x) - \log b, && \forall j\in\mathcal T,\\
& v \ge \sum_{j\in\mathcal T} z_j - k\log b, \qquad v\ge 0, \\
& x \in [1.0001,10^5], \\
& b \in [b_{\mathrm{lb}},2b_{\mathrm{lb}}],
\end{aligned}
\]
\endgroup
where $b_{\mathrm{lb}}=n_{\mathrm{last}}+u_{\mathrm{last}}$, $k=0.15|\mathcal T|$, and $\lambda=10$. Here, the node counts are taken from the latest sampled row. Let $b^\star$ denote the fitted terminal tree size. Using the last $w=5$ sampled points, let $(t_0,n_0)$ and $(t_{\mathrm{last}},n_{\mathrm{last}})$ denote the earliest and latest observations, respectively. The estimator computes the recent time per processed node as $\theta=(t_{\mathrm{last}}-t_0)/(n_{\mathrm{last}}-n_0)$, using a default value of $0.01$ when $n_{\mathrm{last}}\le n_0$. The predicted completion time is then $\widehat T_{\infty}=\tau'+\theta\max\{b^\star-n_{\mathrm{last}},0\}$.

To limit unstable extrapolations, a candidate prediction is accepted only if $\tau'/\widehat T_{\infty}\ge 0.5$. The implementation attempts at most $10$ independently sampled fits and returns the first accepted prediction. If no candidate satisfies the acceptance criterion but finite predictions are available, it returns the smallest predicted value. Every prediction is capped at $10\tau$ to prevent extreme outliers. If all fitting attempts fail, the internal failure value $10^{12}$ is subject to the same cap, assigning the run a completion-time estimate of $10\tau$ ($6{,}000\,\mathrm{s}$ at $\tau=600\,\mathrm{s}$). Such a run is therefore noncompetitive with any solved run. The tuning procedure uses $\widehat T_{\infty}$ to rank configurations, and the following ranking study evaluates its ability to do so. Figure~\ref{fig:svb_flow} summarizes the estimator.

\begin{figure*}[htbp]
  \centering
\caption{Six-step flow of the Case~1 completion-time estimator $\widehat T_{\infty}$.}
  \label{fig:svb_flow}
  \footnotesize
  \begin{tcolorbox}[colback=white,colframe=black!70,boxrule=0.6pt,arc=1.2mm,left=1.8mm,right=1.8mm,top=1.2mm,bottom=1.2mm]
  \centering
  \begin{tabular}{>{\centering\arraybackslash}p{0.93\linewidth}}
  \textbf{Parse the SCIP log and truncate at $\tau'=\min\{\tau,t_{\max}\}$}\\[0.6ex]
  $\Downarrow$\\[0.6ex]
  \textbf{Filter valid rows and discard the startup prefix}\\
  keep only finite rows with time, node count, open nodes, primal bound, and dual bound, after removing rows with \texttt{node}$\le 1$\\[0.6ex]
  $\Downarrow$\\[0.6ex]
  \textbf{Form SVB snapshots}\\
  $g_j=\lvert z_{\mathrm{final}}^{\mathrm{pr}}-z^{\mathrm{du}}(t_j)\rvert$, \quad $\widehat N_j(x)=n_j+u_jx^{g_j}$\\[0.6ex]
  $\Downarrow$\\[0.6ex]
  \textbf{Sample up to 25 rows and fit $(x,b^\star)$ in the log domain}\\
  solve the constrained objective from this subsection on each sampled subset\\[0.6ex]
  $\Downarrow$\\[0.6ex]
  \textbf{Estimate recent time per processed node}\\
  $\theta=(t_{\mathrm{last}}-t_0)/(n_{\mathrm{last}}-n_0)$ over the last $w=5$ sampled points\\[0.6ex]
  $\Downarrow$\\[0.6ex]
  \textbf{Predict completion time and apply safeguards}\\
  $\widehat T_{\infty}=\tau'+\theta\max\{b^\star-n_{\mathrm{last}},0\}$, accept only if $\tau'/\widehat T_{\infty}\ge 0.5$, otherwise retry, fallback to the minimum finite prediction, and cap at $10\tau$
  \end{tabular}
  \end{tcolorbox}
\floatnotes{The estimator is implemented in \texttt{utilities/est\_time.py} and uses only standard solver-log observables. It retries fitting on multiple sampled subsets and assigns the capped value $10\tau$ when all fits fail, making such runs noncompetitive.}
\end{figure*}

\subsection{Retrospective Accuracy on Completed Runs}
Following the offline simulation protocol of \citet{hendel2022estimating}, we construct a set of 100 two-commodity CVRP models with 35 customers, built from instances of the CVRPLIB XML100 collection \citep{lima2014cvrplib,queiroga2021xml100}, each solved to proven optimality within an extended time limit of $1{,}000\,\mathrm{s}$ under the default configuration $p^{\mathrm{def}}$. This yields 100 completed default runs. For each instance with ground-truth time to optimality $T^\star$, we truncate its log at fractions $t_k/T^\star\in\{0.5,0.6,0.7,0.8,0.9\}$ and set $\tau=t_k$. We then apply the estimator using only the information available at the cutoff, including the primal bound at $t_k$, to obtain $\widehat{T}_{\infty}(k)$. We assess prediction accuracy using the ratio $r_k=\widehat{T}_{\infty}(k)/T^\star$ and report the symmetric mean absolute percentage error (SMAPE) and geometric-mean ratio $Q$, where $Q=1$ indicates perfect calibration on the geometric-mean scale. Table~\ref{tab:svb_ratio_dist} summarizes prediction accuracy, the proportion of truncations yielding an accepted fit, and the median fitting time. Prediction error decreases as more of the run is observed, with SMAPE falling from $0.149$ at half of the run to $0.068$ at $90\%$. Accepted fits account for $88\%$--$99\%$ of truncations, and $Q$ lies within $3\%$ of $1$ from the $0.6$ cutoff onward ($0.879$ at the halfway cutoff).

\begin{table*}[htbp]
\centering
\caption{Retrospective accuracy of $\widehat{T}_{\infty}$ on truncated XML100 default-run logs.}
\label{tab:svb_ratio_dist}
\footnotesize
\setlength{\tabcolsep}{5pt}
\renewcommand{\arraystretch}{1.06}
\begin{tabular*}{\textwidth}{@{\extracolsep{\fill}}lrrrrrrr@{}}
\toprule
$t_k/T^\star$ & $n$ & acc. & SMAPE & $Q$ & mean $r$ & med $r$ & med. fit (s) \\
\midrule
0.50 & 100 & 0.88 & 0.149 & 0.879 & 0.889 & 0.909 & 0.104 \\
0.60 & 100 & 0.92 & 0.145 & 0.981 & 0.997 & 1.026 & 0.123 \\
0.70 & 100 & 0.97 & 0.145 & 1.028 & 1.044 & 1.040 & 0.122 \\
0.80 & 100 & 0.99 & 0.110 & 1.022 & 1.035 & 1.009 & 0.115 \\
0.90 & 100 & 0.99 & 0.068 & 1.017 & 1.024 & 1.003 & 0.102 \\
\bottomrule
\end{tabular*}
\floatnotes{Predictions use only information available at each simulated cutoff. The ratio $r=\widehat T_{\infty}/T^\star$ compares predicted and true completion times. ``acc.'' is the proportion of truncations for which a sampled fit satisfies the acceptance rule; the remaining truncations use the minimum-prediction fallback, and no truncation has all fitting attempts fail. SMAPE is $2|\widehat T_{\infty}-T^\star|/(\widehat T_{\infty}+T^\star)$ averaged over all 100 truncations. $Q=\operatorname{GM}(r)$ measures multiplicative calibration, and ``med. fit'' denotes the median fitting time per truncation.}
\end{table*}

\subsection{Ranking Fidelity}
The usefulness of $\widehat{T}_{\infty}$ for tuning depends on its ability to rank configurations by completion time. We evaluate ranking fidelity using XML100 runs collected with a $1{,}000\,\mathrm{s}$ time limit under the default configuration and ten fixed whitelist settings, \texttt{CFG01}--\texttt{CFG10}. Each setting emphasizes a particular solver strategy, such as primal-first search, depth-first search, extensive root cutting, or increased strong branching. The settings files are released with the code. We retain the default configuration and \texttt{CFG02}, \texttt{CFG03}, \texttt{CFG04}, \texttt{CFG05}, \texttt{CFG07}, \texttt{CFG09}, and \texttt{CFG10}. The three more extreme settings, \texttt{CFG01}, \texttt{CFG06}, and \texttt{CFG08}, solve only $26$, $23$, and $24$ of the $100$ instances within $1{,}000\,\mathrm{s}$, respectively. Ground-truth comparisons are restricted to the $49$ instances solved by all eight retained configurations, yielding $392$ configuration runs.

At a simulated cutoff of $600\,\mathrm{s}$, $350$ of these runs have already finished, for which $\widehat{T}_{\infty}$ returns the observed solve times. The remaining $42$ runs require extrapolation. We evaluate ranking on these $42$ runs. For each run, we reconstruct the log summary using only progress rows at or before the cutoff, including the primal bound available at that time, and refit $\widehat{T}_{\infty}$. All $42$ fits succeed. We restrict pairwise comparisons to the $11$ instances with at least two censored configurations, yielding $34$ within-instance configuration pairs. We compare $\widehat{T}_{\infty}$, the final logged gap, and a logged-gap integral, all computed at the same cutoff, against the true completion times. The logged-gap integral serves as a proxy for the primal--dual integral and is computed by left-step integration of the percentage gaps reported in the progress table up to the cutoff, using the first visible gap for the initial interval and $100$ for missing gaps. This proxy differs from the primal--dual integral reported by \textsc{SCIP} and used in the Case~2 loss. Table~\ref{tab:rank_corr} reports the number of correctly ordered pairs and the mean within-instance Kendall rank correlation \citep{kendall1938new,kendall1945treatment}
\[
\tau_b=\frac{C-D}{\sqrt{(C+D+T_x)(C+D+T_y)}},
\]
where $C$ and $D$ denote the numbers of concordant and discordant pairs, respectively, and $T_x,T_y$ count ties in the two rankings.

\begin{table}[htbp]
\centering
\caption{Prediction-only ranking at a $600\,\mathrm{s}$ cutoff.}
\label{tab:rank_corr}
\footnotesize
\setlength{\tabcolsep}{6pt}
\renewcommand{\arraystretch}{1.05}
\begin{tabular}{lcc}
\toprule
\textnormal{Signal at the cutoff} & \textnormal{Correctly ordered pairs} & \textnormal{Mean within-instance $\tau_b$} \\
\midrule
Completion estimate $\widehat{T}_{\infty}$ & \textbf{26/34} & \textbf{0.515} \\
Final logged gap & 25/34 & 0.485 \\
Logged-gap integral (PDI proxy) & 17/34 & 0.121 \\
\bottomrule
\end{tabular}
\floatnotes{The evaluation includes the 11 XML100 instances with at least two censored configurations. All ranking signals are computed using only information available at the cutoff and are compared on the same 34 within-instance configuration pairs.}
\end{table}

The completion-time estimate places unfinished runs on the same time scale as the observed solve times of completed runs, so a single Case~1 objective ranks solved and unsolved configurations together, which neither gap-based signal provides. Among the three signals, it also orders the most configuration pairs correctly ($26$ of $34$, against $17$ for the logged-gap integral) and achieves the highest mean within-instance rank correlation.

\section{Additional Results, Transfer Breakdowns, and Hard-Regime Audit}\label{app_sensitivity}
This section provides supplementary analyses for the experiments in Section~\ref{sec_exp}, including the edit sizes of Exp.~A, split-level and benchmark-level breakdowns of Exp.~C, and an audit of the Case~3 outlier in Exp.~A.

\paragraph{Exp.~A edit sizes.}
Table~\ref{tab:edit_size} reports the size of the edit sets that \name{} proposes in Exp.~A, computed from the archived tuning histories of all $215$ instances. Each proposal is measured against the best configuration at the time it is made. Edits are smallest in Case~1 and largest in Case~3.

\begin{table}[htbp]
\centering
\caption{Edit-set sizes of \name{} proposals in Exp.~A.}
\label{tab:edit_size}
\footnotesize
\setlength{\tabcolsep}{4pt}
\renewcommand{\arraystretch}{1.05}
\begin{tabular}{lrrrrrrr}
\toprule
\textnormal{Case} & \textnormal{Proposals} & \textnormal{Median} & \textnormal{Mean} & \textnormal{90th pct.} & \textnormal{Max} & \textnormal{$\le 10$ (\%)} & \textnormal{At cap (\%)} \\
\midrule
Case~1 & $3{,}800$ & 5 & 5.00 & 8 & 19 & 97.8 & 0.0 \\
Case~2 & $4{,}550$ & 7 & 7.35 & 12 & 20 & 86.5 & 0.4 \\
Case~3 & $2{,}400$ & 11 & 11.71 & 18 & 20 & 49.3 & 4.7 \\
\midrule
All & $10{,}750$ & 7 & 7.49 & 13 & 20 & 82.2 & 1.2 \\
\bottomrule
\end{tabular}
\floatnotes{Edit-set size is the number of whitelisted parameters for which the evaluated configuration differs from the best configuration at the time of the proposal. Each proposal may edit at most $20$ of the $45$ whitelisted parameters.}
\end{table}

\paragraph{Exp.~C split-level and benchmark-level breakdowns.}
The aggregate results of Exp.~C conceal substantial variation across benchmark families. For Case~1 test loss, joint tuning outperforms independent selection at $m=2$ and $m=3$, achieves nearly identical performance at $m=4$, and performs worse on the eligible pooled subsets at $m=5$ and $m=6$. Table~\ref{tab:app_gen_case1_by_benchmark} reports the Case~1 test loss by benchmark family and training size, while Table~\ref{tab:app_gen_transfer_by_benchmark} provides the transfer results across regimes at $m=3$. Joint tuning consistently outperforms independent selection on CVRP across all five training sizes and on GAP at both eligible sizes. In contrast, independent selection performs better on UFL at every training size. Results for JSP are mixed, with independent selection slightly better at $m=2$ and joint tuning better at $m=3$. The pooled results at larger $m$ therefore reflect the contrasting performance of the two remaining eligible families, CVRP and UFL.

\begin{table*}[htbp]
\centering
\caption{Benchmark-specific Case~1 test losses across training sizes in Exp.~C.}
\label{tab:app_gen_case1_by_benchmark}
\footnotesize
\setlength{\tabcolsep}{3pt}
\renewcommand{\arraystretch}{1.05}
\begin{tabular*}{\textwidth}{@{\extracolsep{\fill}}l cc cc cc cc cc@{}}
\toprule
\textnormal{Benchmark} & \multicolumn{2}{c}{$m=2$} & \multicolumn{2}{c}{$m=3$} & \multicolumn{2}{c}{$m=4$} & \multicolumn{2}{c}{$m=5$} & \multicolumn{2}{c}{$m=6$} \\
\cmidrule(lr){2-3}\cmidrule(lr){4-5}\cmidrule(lr){6-7}\cmidrule(lr){8-9}\cmidrule(l){10-11}
& Indep. & Joint & Indep. & Joint & Indep. & Joint & Indep. & Joint & Indep. & Joint \\
\midrule
CVRP & 1.122 & \textbf{0.796} & 1.014 & \textbf{0.874} & 0.985 & \textbf{0.885} & 1.028 & \textbf{0.877} & 0.769 & \textbf{0.675} \\
GAP  & 0.775 & \textbf{0.673} & 0.841 & \textbf{0.680} & -- & -- & -- & -- & -- & -- \\
JSP  & \textbf{0.265} & 0.325 & 0.217 & \textbf{0.204} & -- & -- & -- & -- & -- & -- \\
UFL  & \textbf{0.369} & 0.382 & \textbf{0.323} & 0.395 & \textbf{0.318} & 0.355 & \textbf{0.287} & 0.445 & \textbf{0.292} & 0.382 \\
\bottomrule
\end{tabular*}
\floatnotes{Each entry is the geometric mean of split-level Case~1 test losses for the corresponding benchmark family, tuning strategy, and training size. A dash (--) indicates that the family is ineligible for the given training size under the requirement of at least four test instances. Lower values indicate better performance.}
\end{table*}

\begin{table*}[htbp]
\centering
\caption{Benchmark-specific transfer results in Exp.~C at $m=3$.}
\label{tab:app_gen_transfer_by_benchmark}
\footnotesize
\setlength{\tabcolsep}{5pt}
\renewcommand{\arraystretch}{1.05}
\begin{tabular*}{\textwidth}{@{\extracolsep{\fill}}lcccccc@{}}
\toprule
\textnormal{Benchmark} & \textnormal{Indep. C1} & \textnormal{Joint C1} & \textnormal{Indep. C2} & \textnormal{Joint C2} & \textnormal{Indep. C3} & \textnormal{Joint C3} \\
\midrule
CVRP & 1.014 & \textbf{0.874} & 0.776 & \textbf{0.749} & 1.018 & \textbf{1.011} \\
GAP  & 0.841 & \textbf{0.680} & 1.030 & \textbf{0.913} & -- & -- \\
JSP  & 0.217 & \textbf{0.204} & \textbf{0.807} & 0.830 & -- & -- \\
UFL  & \textbf{0.323} & 0.395 & -- & -- & -- & -- \\
\bottomrule
\end{tabular*}
\floatnotes{Columns C1, C2, and C3 report Case~1 test loss, same-family Case~2 loss, and same-family Case~3 loss, respectively. Each entry is the geometric mean over repeated splits for the corresponding benchmark family and tuning strategy. A dash (--) indicates that the family has no test instances in that regime. Lower values indicate better performance.}
\end{table*}

\paragraph{Case~3 outlier audit.}
The aggregate Case~3 results in Exp.~A are strongly influenced by a single instance, \texttt{miplib\_irish-electricity}. Table~\ref{tab:case3_outliers} reports the ten instances with the largest absolute differences in the selected Case~3 loss between \name{} and \textsc{RBFOpt}, using the same best-so-far selection rule as Table~\ref{tab:exp_main_results}. The largest difference occurs on \texttt{miplib\_irish-electricity}, where \textsc{RBFOpt} achieves a loss lower by $0.256$. On \texttt{miplib\_neos-1354092}, \name{} achieves a loss lower by $0.079$, while all remaining differences are below $0.007$. Across all $48$ instances, the geometric-mean losses are $0.9931$, $0.9931$, $0.9882$, and $0.9929$ for \name{}, \textsc{SMAC3}, \textsc{RBFOpt}, and \textsc{Random}, respectively. Excluding only \texttt{miplib\_irish-electricity} yields losses of $0.9939$, $0.9947$, $0.9955$, and $0.9950$, respectively. The aggregate advantage of \textsc{RBFOpt} in Case~3 therefore comes from its improvement on a single instance. The main results in Table~\ref{tab:exp_main_results} retain all $48$ instances.

\begin{table}[htbp]
\centering
\caption{Largest instance-level loss differences between \name{} and \textsc{RBFOpt} in Case~3.}
\label{tab:case3_outliers}
\footnotesize
\setlength{\tabcolsep}{4pt}
\renewcommand{\arraystretch}{1.05}
\begin{tabular}{L{0.32\linewidth}ccccc}
\toprule
\textnormal{Instance} & \textnormal{\name{}} & \textnormal{SMAC3} & \textnormal{RBFOpt} & \textnormal{Random} & $\Delta$ \\
\midrule
\texttt{miplib\_irish-electricity} & 0.9533 & 0.9218 & 0.6977 & 0.8991 & $+0.2556$ \\
\texttt{miplib\_neos-1354092} & 0.9215 & 0.9863 & 1.0000 & 1.0000 & $-0.0785$ \\
\texttt{cvrp\_B-n51-k7} & 0.9922 & 0.9804 & 0.9855 & 0.9860 & $+0.0067$ \\
\texttt{cvrp\_B-n45-k6} & 0.9942 & 0.9893 & 0.9894 & 0.9895 & $+0.0048$ \\
\texttt{cvrp\_A-n62-k8} & 0.9964 & 0.9918 & 0.9919 & 0.9918 & $+0.0046$ \\
\texttt{cvrp\_P-n50-k10} & 0.9919 & 0.9959 & 0.9960 & 0.9960 & $-0.0041$ \\
\texttt{cvrp\_B-n67-k10} & 0.9964 & 0.9920 & 0.9923 & 0.9919 & $+0.0041$ \\
\texttt{cvrp\_P-n65-k10} & 0.9927 & 0.9958 & 0.9967 & 0.9954 & $-0.0041$ \\
\texttt{cvrp\_A-n63-k9} & 0.9953 & 0.9993 & 0.9993 & 0.9993 & $-0.0040$ \\
\texttt{cvrp\_A-n69-k9} & 0.9972 & 0.9941 & 0.9932 & 0.9936 & $+0.0040$ \\
\bottomrule
\end{tabular}
\floatnotes{The comparison covers the 48 Case~3 instances, with the default configuration eligible for selection by every method. Differences are defined as $\Delta=L_{\name{}}-L_{\textsc{RBFOpt}}$, so positive $\Delta$ indicates a lower loss for \textsc{RBFOpt}. All reported losses are the original Case~3 values without display magnification.}
\end{table}

\section{Code and Data Availability, Environment, Compute Accounting, and Licensing}\label{app_artifact}
\paragraph{Code and data.}
The \name{} source code is available at \coderepo. The repository includes the tuning loop, typed parameter whitelist, prompt construction, log parser, completion-time estimator, baseline implementations, experiment entry points, benchmark manifests, split assignments, representative prompts, and analysis scripts. It also provides instructions for obtaining the original benchmark instances and regenerating and verifying the derived \texttt{.mps} files, along with regime assignments and file-level provenance. The repository also includes per-instance results for Exp.~A, Exp.~B, C1, and the XML100 validation studies, together with a script that verifies the corresponding reported summaries.

\paragraph{Licenses and upstream assets.}
The benchmark instances are obtained from the public sources cited in Appendix~\ref{app_instance} and remain subject to their respective upstream terms. For each source, the repository records the citation, upstream URL, and licensing status, including cases where the publicly available files do not include a separate license. The generated \texttt{.mps} files are deterministic reformulations of the original benchmark instances, and their source collections are documented in the repository. The \name{} source code is released under the MIT license. Table~\ref{tab:artifact_env} summarizes the hardware, software, and token-side compute accounting used for reproducibility.

\begin{table}[htbp]
\centering
\caption{Hardware, software, and LLM token accounting for the reported experiments.}
\label{tab:artifact_env}
\footnotesize
\setlength{\tabcolsep}{4pt}
\renewcommand{\arraystretch}{1.05}
\begin{tabular}{@{}p{0.43\linewidth}p{0.49\linewidth}@{}}
\toprule
Item & Value \\
\midrule
Operating system & Linux (Ubuntu; kernel 6.14.0-33-generic; x86\_64) \\
CPU & 13th Gen Intel(R) Core(TM) i7-13700K (24 logical CPUs) \\
System memory & 62 GiB total RAM \\
Python & 3.12.4 \\
SCIP & 9.2.4 \\
Threads per solver run & 1 \\
Solver-side seed convention & \texttt{randomization/randomseedshift = 42} \\
LLM model & GPT-5.2 \\
Total LLM calls in Exp.~A & 10,964 \\
Total tokens in Exp.~A & 62,958,576 \\
Approx. uncached list-price cost & USD 320.08 \\
Invalid-JSON warnings in Exp.~A & 69 \\
\bottomrule
\end{tabular}
\end{table}

\paragraph{Compute accounting convention.}
The main text reports the exact solver-run budgets, and exact token accounting for Exp.~A is reported in Table~\ref{tab:artifact_env}. All solver calls are single-threaded and blocking. Exp.~A uses 51 solver evaluations per method-instance pair, comprising one default run and $K=50$ proposals. Exp.~B uses 31 evaluations per method-instance pair, comprising one default run and $K=30$ proposals. Exp.~C uses the same budget per training run for independent and joint tuning, with additional evaluations on test instances determined by the fixed split manifests.

\paragraph{Seed convention.}
All reported solver runs fix \texttt{randomization/randomseedshift} to 42, and the stochastic black-box baselines use an optimizer seed of 42. This convention makes every solver run reproducible. The LLM policy uses a fixed model and API parameter settings, as documented in Appendix~\ref{app_hparams}.

\section{Full Benchmark Manifests and Family-Specific Formulations}\label{app_instance_lists}
This appendix provides the family-specific formulations, source breakdowns, and benchmark manifests supporting Appendix~\ref{app_instance}. The experimental benchmark comprises 90 CVRP, 12 GAP, 24 JSP, 30 UFL, and 59 MIPLIB instances, 215 in total. Appendix~\ref{app_instance} describes the benchmark sources, construction pipeline, regime assignment, and split protocol; this appendix provides the detailed mathematical formulations and source tables for each benchmark family.

\subsection{CVRP: Capacitated Vehicle Routing Problem}\label{app_sub_cvrp}
\paragraph{Problem definition.}
The CVRP seeks a set of vehicle routes of minimum total cost, each starting and ending at a depot, such that every customer is served exactly once and each vehicle respects the capacity limit $Q$. Let $\mathcal{G}=(\mathcal{V},\mathcal{E})$ be a complete undirected graph with two copies of the depot, where $\mathcal{V}=\{0,1,\ldots,n,n{+}1\}$ and $\mathcal{C}=\mathcal{V}\setminus\{0,n{+}1\}$ denotes the customer set. Each customer $i$ has demand $d_i$, and $K$ vehicles are available.

\paragraph{Formulation.}
Let $c_{ij}$ denote the cost of using edge $\{i,j\}$. We adopt the two-commodity flow formulation of \citet{baldacci2004exact}, with binary edge variables $x_{ij}$ and nonnegative continuous flow variables $y_{ij}\ge 0$:
\begin{equation*}
\begin{alignedat}{2}
\min\quad & \sum_{\{i,j\}\in \mathcal{E}} c_{ij}\,x_{ij} & & \\
\text{s.t.}\quad
& \sum_{j\in \mathcal{V}\setminus\{i\}} \big(y_{ji} - y_{ij}\big) = 2d_i,
&\quad& \forall i\in \mathcal{C}, \\
& \sum_{j\in \mathcal{C}} y_{0j} = d(\mathcal{C}), & & \\
& \sum_{j\in \mathcal{C}} y_{j0} = KQ - d(\mathcal{C}), & & \\
& \sum_{j\in \mathcal{C}} y_{n+1,j} = KQ, & & \\
& y_{ij} + y_{ji} = Q\,x_{ij},
&\quad& \forall \{i,j\}\in \mathcal{E}, \\
& \sum_{\substack{j\in \mathcal{V}\\ j>i}} x_{ij}
 + \sum_{\substack{j\in \mathcal{V}\\ j<i}} x_{ji} = 2,
&\quad& \forall i\in \mathcal{C}, \\
& y_{ij} \ge 0,\; y_{ji} \ge 0,
&\quad& \forall \{i,j\}\in \mathcal{E}, \\
& x_{ij} \in \{0,1\},
&\quad& \forall \{i,j\}\in \mathcal{E}.
\end{alignedat}
\end{equation*}

Here $d(\mathcal{C})=\sum_{i\in\mathcal{C}} d_i$. We use the 90 CVRPLIB instances summarized in Table~\ref{tab:cvrplib_sets_summary}.

\paragraph{Instance selection.}
We select CVRP instances from the \texttt{A}, \texttt{B}, \texttt{E}, \texttt{F}, and \texttt{P} families in CVRPLIB \citep{lima2014cvrplib}, originating from \citet{augerat1995computational} (sets A/B/P), \citet{christofides1969algorithm} (set E), and \citet{fisher1994optimal} (set F). CVRPLIB instance names follow the format \texttt{<set>-n<}$N$\texttt{>-k<}$K$\texttt{>}, where $N$ denotes the number of nodes, including the depot, and $K$ is the number of routes in the reference solution. For an instance with $N$ nodes, our implementation sets $n:=N-1$ customers, represents the depot by nodes $0$ and $n{+}1$, and fixes the fleet size to $K$. Table~\ref{tab:cvrplib_sets_summary} summarizes the family counts and size ranges of the 90 selected instances.
\begin{table}[htbp]
\centering
\caption{Summary of CVRP instances used in the experiments.}
\label{tab:cvrplib_sets_summary}
\footnotesize
\setlength{\tabcolsep}{5pt}
\renewcommand{\arraystretch}{1.05}
\begin{tabular*}{\linewidth}{@{\extracolsep{\fill}}lccc@{}}
\toprule
\textnormal{Set} & \textnormal{\#} & \textnormal{Range of $N$} & \textnormal{Source} \\
\midrule
\texttt{A} & 27 & 32--80   & \citet{augerat1995computational} \\
\texttt{B} & 23 & 31--78   & \citet{augerat1995computational} \\
\texttt{E} & 13 & 13--101  & \citet{christofides1969algorithm} \\
\texttt{F} & 3  & 45--135  & \citet{fisher1994optimal} \\
\texttt{P} & 24 & 16--101  & \citet{augerat1995computational} \\
\midrule
\textbf{Total} & \textbf{90} & & \\
\bottomrule
\end{tabular*}
\end{table}

\subsection{GAP: Generalized Assignment Problem}\label{app_sub_gap}
\paragraph{Problem definition.}
The GAP assigns each job $j\in\mathcal{J}$ to exactly one agent $i\in\mathcal{I}$. Assigning job $j$ to agent $i$ consumes $a_{ij}$ units of the agent's capacity $b_i$ and incurs cost $c_{ij}$. The objective is to minimize the total assignment cost while respecting all agent capacities \citep{ross1975branch}. The binary variable $x_{ij}$ indicates whether job $j$ is assigned to agent $i$.

\paragraph{Formulation.}
\begin{equation*}
\begin{alignedat}{2}
\min\quad & \sum_{i\in \mathcal{I}}\sum_{j\in \mathcal{J}} c_{ij}\,x_{ij} & & \\
\text{s.t.}\quad
& \sum_{i\in \mathcal{I}} x_{ij} = 1, &\quad& \forall j\in \mathcal{J},\\
& \sum_{j\in \mathcal{J}} a_{ij}\,x_{ij} \le b_i, &\quad& \forall i\in \mathcal{I},\\
& x_{ij}\in\{0,1\}, &\quad& \forall i\in \mathcal{I},\ \forall j\in \mathcal{J}.
\end{alignedat}
\end{equation*}

We use 12 OR-Library instances from \texttt{gapc} and \texttt{gapd}.

\paragraph{Instance selection.}
We select the large-scale type-C and type-D instances from the OR-Library GAP collection \citep{beasley1990library}, introduced by \citet{chu1997genetic} and distributed as \texttt{gapc} and \texttt{gapd}. Each file contains six minimization instances with $|\mathcal{I}|\in\{5,10,20\}$ and $|\mathcal{J}|\in\{100,200\}$. Instances are identified by their source file and within-file index (e.g., \texttt{gapc-3}). Table~\ref{tab:gap_summary} summarizes the two source files and their instance sizes.
\begin{table}[htbp]
\centering
\caption{Summary of GAP instances used in the experiments.}
\label{tab:gap_summary}
\footnotesize
\setlength{\tabcolsep}{5pt}
\renewcommand{\arraystretch}{1.05}
\begin{tabular*}{\linewidth}{@{\extracolsep{\fill}}lccc@{}}
\toprule
\textnormal{File} & \textnormal{\#} & \textnormal{Sizes $(|\mathcal{I}|\times |\mathcal{J}|)$} & \textnormal{Source} \\
\midrule
\texttt{gapc}  & 6 & $\{5,10,20\}\times\{100,200\}$ & \multirow{2}{*}{\citet{chu1997genetic}} \\
\texttt{gapd}  & 6 & $\{5,10,20\}\times\{100,200\}$ & \\
\midrule
\textbf{Total} & \textbf{12} & & \\
\bottomrule
\end{tabular*}
\end{table}

\subsection{JSP: Job Shop Scheduling Problem}\label{app_sub_jsp}
\paragraph{Problem definition.}
Let $\mathcal J$ and $\mathcal M$ denote the sets of jobs and machines. The JSP schedules jobs, each comprising an ordered sequence of operations. Each operation $o$ requires processing time $p_o$ on a designated machine $\mu(o)$. Each machine can process at most one operation at a time, and preemption is not allowed. The objective is to minimize the makespan $C_{\max}$ \citep{manne1960job}.

\paragraph{Formulation.}
We use the valid scheduling horizon $H=\sum_{o} p_o$. Let $\mathcal{P}$ denote the within-job precedence arcs and $\mathcal{D}$ the set of unordered pairs of operations requiring the same machine. Using start-time variables $s_o$ and binary ordering variables $\delta_{oo'}$ for $(o,o')\in\mathcal{D}$, we adopt the following standard disjunctive MILP formulation:
\begin{equation*}
\begin{alignedat}{2}
\min\quad & C_{\max} & & \\
\text{s.t.}\quad
& s_{o'} \ge s_o + p_o, &\quad& \forall (o,o')\in\mathcal{P},\\
& s_o+p_o \le s_{o'} + H\,(1-\delta_{oo'}), &\quad& \forall (o,o')\in\mathcal{D},\\
& s_{o'}+p_{o'} \le s_o + H\,\delta_{oo'}, &\quad& \forall (o,o')\in\mathcal{D},\\
& C_{\max} \ge s_o + p_o, &\quad& \forall o,\\
& s_o\ge 0, &\quad& \forall o,\\
& \delta_{oo'}\in\{0,1\}, &\quad& \forall (o,o')\in\mathcal{D}.
\end{alignedat}
\end{equation*}

We use 24 instances from the OR-Library \texttt{jobshop1} collection, covering the \texttt{abz}, \texttt{ft}, \texttt{la}, and \texttt{orb} families.

\paragraph{Instance selection.}
We select 24 benchmark instances from the OR-Library \texttt{jobshop1} collection \citep{beasley1990library}, spanning the \texttt{abz}, \texttt{ft}, \texttt{la}, and \texttt{orb} families \citep{adams1988shifting,fisher1963probabilistic,lawrence1984supplement,applegate1991computational}. Table~\ref{tab:jsp_sets_summary} summarizes the number of instances, size ranges, and sources for each family.
\begin{table}[htbp]
\centering
\caption{Summary of JSP instances used in the experiments.}
\label{tab:jsp_sets_summary}
\footnotesize
\setlength{\tabcolsep}{5pt}
\renewcommand{\arraystretch}{1.05}
\begin{tabular*}{\linewidth}{@{\extracolsep{\fill}}lccc@{}}
\toprule
\textnormal{Family} & \textnormal{\#} & \textnormal{Range of $(|\mathcal{J}|,|\mathcal{M}|)$} & \textnormal{Source} \\
\midrule
\texttt{abz(5-9)}   & 5  & $(10$--$20,\ 10$--$15)$ & \citet{adams1988shifting} \\
\texttt{ft(10,20)}  & 2  & $(10$--$20,\ 5$--$10)$  & \citet{fisher1963probabilistic} \\
\texttt{la(01-12)}  & 12 & $(10$--$20,\ 5)$        & \citet{lawrence1984supplement} \\
\texttt{orb(01-05)} & 5  & $(10,\ 10)$             & \citet{applegate1991computational} \\
\midrule
\textbf{Total} & \textbf{24} & & \\
\bottomrule
\end{tabular*}
\end{table}

\subsection{UFL: Uncapacitated Facility Location Problem}\label{app_sub_facility}
\paragraph{Problem definition.}
The UFL seeks to minimize the total facility opening and customer assignment costs, where each facility incurs a fixed opening cost and every customer must be assigned to an open facility \citep{kuehn1963heuristic}. Let $\mathcal F$ and $\mathcal C$ denote the facility and customer sets, respectively. Opening facility $j$ incurs cost $f_j$, and assigning customer $i$ to facility $j$ incurs cost $c_{ij}$. The binary variables $y_j$ and $x_{ij}$ indicate whether facility $j$ is open and whether customer $i$ is assigned to it, respectively.

\paragraph{Formulation.}
\begin{equation*}
\begin{alignedat}{2}
\min\quad & \sum_{j\in \mathcal{F}} f_j\,y_j \;+\; \sum_{i\in \mathcal{C}}\sum_{j\in \mathcal{F}} c_{ij}\,x_{ij} & & \\
\text{s.t.}\quad
& \sum_{j\in \mathcal{F}} x_{ij} = 1, &\quad& \forall i\in \mathcal{C},\\
& x_{ij} \le y_j, &\quad& \forall i\in \mathcal{C},\ \forall j\in \mathcal{F},\\
& x_{ij}\in\{0,1\}, &\quad& \forall i\in \mathcal{C},\ \forall j\in \mathcal{F},\\
& y_j\in\{0,1\}, &\quad& \forall j\in \mathcal{F}.
\end{alignedat}
\end{equation*}

We use 30 UflLib instances from the Uniform family.

\paragraph{Instance selection.}
We select uncapacitated facility location instances from UflLib \citep{hoefer2018ufllib}. Following our benchmark protocol, we include the complete Uniform family, comprising all 30 \texttt{*UnifS} instances. We exclude the \texttt{*EuclS} instances, which are typically solved at the root node under default settings and therefore provide limited information for tuning. Table~\ref{tab:ufl_sets_summary} summarizes the UFL instances included in the benchmark.
\begin{table}[htbp]
\centering
\caption{Summary of UFL instances used in the experiments.}
\label{tab:ufl_sets_summary}
\footnotesize
\setlength{\tabcolsep}{5pt}
\renewcommand{\arraystretch}{1.05}
\begin{tabular*}{\linewidth}{@{\extracolsep{\fill}}lccc@{}}
\toprule
\textnormal{Family} & \textnormal{\#} & \textnormal{$(|\mathcal{F}|,|\mathcal{C}|)$} & \textnormal{Source} \\
\midrule
\texttt{*UnifS} & 30 & $(100,\ 100)$ & \citet{hoefer2018ufllib} \\
\midrule
\textbf{Total} & \textbf{30} & & \\
\bottomrule
\end{tabular*}
\end{table}

\subsection{General MILP Benchmark: MIPLIB 2017}\label{app_sub_miplib}
\paragraph{Benchmark description.}
To complement the domain-specific models, we include instances from the official MIPLIB~2017 benchmark set \citep{gleixner2021miplib}. MIPLIB~2017 is the sixth release of the library and was constructed through a data-driven selection process. From an initial pool of 5{,}721 submissions, it provides a collection of 1{,}065 instances, including a curated benchmark subset of 240 instances for solver evaluation (\texttt{benchmark-v2}).

\paragraph{Instance selection.}
We uniformly sample $60$ instances from the official MIPLIB~2017 benchmark subset using a fixed random seed ($42$). We exclude \texttt{k1mushroom} because it does not complete presolving under default SCIP settings within the 600-second time limit. The final experimental subset therefore contains $59$ instances, listed in Table~\ref{tab:miplib_subset}.

\begin{table}[htbp]
  \centering
\caption{Final 59-instance MIPLIB~2017 benchmark subset.}
\label{tab:miplib_subset}
  \scriptsize
  \setlength{\tabcolsep}{5pt}
  \renewcommand{\arraystretch}{1.05}
    \begin{tabular*}{\textwidth}{@{\extracolsep{\fill}}lll@{}}
      \toprule
      \multicolumn{3}{l}{\textnormal{Final selected instances (\texttt{benchmark-v2})}} \\
      \midrule
      \texttt{50v-10} & \texttt{assign1-5-8} & \texttt{atlanta-ip} \\
      \texttt{b1c1s1} & \texttt{beasleyC3} & \texttt{cbs-cta} \\
      \texttt{chromaticindex512-7} & \texttt{cmflsp50-24-8-8} & \texttt{CMS750\_4} \\
      \texttt{cod105} & \texttt{comp21-2idx} & \texttt{cryptanalysiskb128n5obj16} \\
      \texttt{dano3\_3} & \texttt{drayage-25-23} & \texttt{dws008-01} \\
      \texttt{fhnw-binpack4-48} & \texttt{gfd-schedulen180f7d50m30k18} & \texttt{glass4} \\
      \texttt{glass-sc} & \texttt{gmu-35-50} & \texttt{h80x6320d} \\
      \texttt{irish-electricity} & \texttt{lectsched-5-obj} & \texttt{momentum1} \\
      \texttt{mushroom-best} & \texttt{mzzv11} & \texttt{n3div36} \\
      \texttt{neos8} & \texttt{neos859080} & \texttt{neos-662469} \\
      \texttt{neos-827175} & \texttt{neos-911970} & \texttt{neos-1354092} \\
      \texttt{neos-1445765} & \texttt{neos-3046615-murg} & \texttt{neos-3083819-nubu} \\
      \texttt{neos-3627168-kasai} & \texttt{neos-3988577-wolgan} & \texttt{neos-5052403-cygnet} \\
      \texttt{nexp-150-20-8-5} & \texttt{nursesched-sprint02} & \texttt{p200x1188c} \\
      \texttt{piperout-27} & \texttt{radiationm18-12-05} & \texttt{radiationm40-10-02} \\
      \texttt{ran14x18-disj-8} & \texttt{rocI-4-11} & \texttt{rocII-5-11} \\
      \texttt{roll3000} & \texttt{s100} & \texttt{snp-02-004-104} \\
      \texttt{supportcase22} & \texttt{supportcase40} & \texttt{supportcase42} \\
      \texttt{swath1} & \texttt{thor50dday} & \texttt{traininstance2} \\
      \texttt{uct-subprob} & \texttt{unitcal\_7} & \\
      \bottomrule
    \end{tabular*}
\floatnotes{The subset used in the experiments was sampled from the official \texttt{benchmark-v2} collection after removing \texttt{k1mushroom}, which does not finish presolving under the default 600-second run.}
\end{table}

\end{document}